%% file: mainPaper.tex
\documentclass[final,leqno,onefignum,onetabnum]{siamltexmm}

\usepackage{amssymb}
\usepackage{tikz}
\usepackage{stackrel}
\usepackage{url}
\usepackage{xspace}
\usepackage{multirow}
\usepackage[T1]{fontenc}
\usepackage{lmodern}
\usepackage{amsmath}
\usepackage{amsfonts}
\usepackage{graphicx}
\usepackage{bm}
\usepackage{subfigure}
\usepackage{comment}

\usepackage[letterpaper,centering]{geometry}
\usepackage[textsize=tiny]{todonotes}

\usepackage{ifpdf}
\ifpdf
  \usepackage{epstopdf} 
\fi

\newcommand{\mat}[1]{\mathbf{#1}}
\newcommand{\pmat}[2]{\bm{\Phi}_{#1:#2}}
\newcommand{\pvec}[1]{\bm{\phi}_{#1}}
\newcommand{\psvec}[1]{\bm{\psi}_{#1}}
\newcommand{\mbf}[1]{\bm{#1}}

\newcommand{\varn}[0]{\sigma^2} 

\DeclareMathOperator*{\Span}{span}

\newcommand{\sedit}{\color{black}}
\newcommand{\dedit}{\color{black}}

\newcommand{\ssedit}{\color{black}}
\newcommand{\ddedit}{\color{black}}

\title{Mercer kernels and integrated variance experimental design: connections between Gaussian process regression and polynomial approximation}
\author{Alex Gorodetsky\footnotemark[1] \and Youssef Marzouk\footnotemark[1]}

\begin{document}

\maketitle

\renewcommand{\thefootnote}{\fnsymbol{footnote}}

\footnotetext[1]{Department of Aeronautics and Astronautics, Massachusetts Institute of Technology, Cambridge, MA
02139, USA. \texttt{\{goroda,ymarz\}@mit.edu}.}

\renewcommand{\thefootnote}{\arabic{footnote}}

\newcommand{\slugmaster}{%
\slugger{xxx}{xxxx}{xx}{x}}

\begin{abstract}
This paper examines experimental design procedures used to develop surrogates of computational models, exploring the interplay between experimental designs and approximation algorithms. We focus on two widely used approximation approaches, Gaussian process (GP) regression and non-intrusive polynomial approximation. 
First, we introduce algorithms for minimizing a posterior {integrated variance} (IVAR) design criterion for GP regression. Our formulation treats design as a continuous optimization problem that can be solved {with gradient-based methods} on complex input domains, without resorting to greedy approximations. We show that minimizing IVAR in this way yields point sets with good interpolation properties, and that it enables more accurate GP regression than designs based on entropy minimization or mutual information maximization. 
Second, using a Mercer kernel/eigenfunction perspective on GP regression, we identify conditions under which GP regression coincides with pseudospectral polynomial approximation. Departures from these conditions can be understood as changes either to the kernel or to the experimental design itself. We then show how IVAR-optimal designs, while sacrificing discrete orthogonality of the kernel eigenfunctions, can yield lower approximation error than
orthogonalizing point sets. 
Finally, we compare the performance of adaptive Gaussian process regression and adaptive pseudospectral approximation for several classes of target functions, identifying features that are favorable to the GP + IVAR approach.

\end{abstract}


\pagestyle{myheadings}
\thispagestyle{plain}
\markboth{\MakeUppercase{A.\ Gorodetsky and Y.\ M.\ Marzouk}}{\MakeUppercase{IVAR Experimental Design}}

\input{intro2}

\input{gp}

\input{expDesign2}

\input{numExamples}

\input{psa}

\input{theory}
\input{GPPSAexp}
\input{conclusion}

\subsection*{Acknowledgments}
The authors gratefully acknowledge BP for funding this research. We
thank Tarek Moselhy, Akil Narayan, and Alessio Spantini for many
helpful discussions. We also thank the referees and editors for their
thoughtful comments and suggestions.

\begin{appendix}
\input{proofsContinuousVersion}

\input{gradient_appendix}
\end{appendix}

\bibliographystyle{plain}
\bibliography{main}

\end{document}

%% file: intro2.tex
\section{Introduction}

Computational simulations are essential for design, optimization, uncertainty quantification, and inference in complex systems. Yet these tasks typically require a large number of simulations over a range of parameter values, which can be computationally prohibitive. One method of mitigating this computational expense is to construct surrogates or ``emulators'' that replace the simulation in the relevant analyses. Because many computational simulations are available only as black-box or legacy codes, surrogates often must be constructed through a limited number of simulations at particular parameter values $\{x_i\}$. These simulations are sometimes called ``computer experiments'' \cite{SacksWelch1989, Santner2003} and choosing these parameter values is a question of experimental design. 

Surrogate construction can be viewed as a function approximation problem in which one attempts to approximate an input-output relationship $f(x)$ induced by the expensive simulation. One can do so deterministically, i.e., obtaining a single approximation $\hat{f}$, or probabilistically, i.e., obtaining a distribution over possible functions $\hat{f} \sim \mathcal{F}$, where $\mathcal{F}$ is a probability distribution on a suitable function space. In either case, three decisions must be made. First, one must choose an \textit{approximation space} that contains candidate surrogate functions; second, one must select a set of parameter values or \textit{experiments} at which to simulate the system; and finally, one must choose an \textit{algorithm} to convert the simulation results into a particular function $\hat{f}$ or distribution $\mathcal{F}$. 
Examples of approximation spaces include those spanned by polynomials of a certain degree, or by radial basis functions and other kernels. Possible experiments include designs produced by Monte Carlo sampling, obtained via optimization of information-based design criteria such as entropy and mutual information, or based on numerical quadrature rules. Finally, examples of algorithms include linear regression and pseudospectral approximation. These three decisions are usually not independent; for example, pseudospectral approximation requires experimental design procedures that preserve orthogonality of the relevant basis functions.

\sedit
In this paper, we analyze two of the most commonly used surrogate construction approaches in uncertainty quantification: Gaussian process regression (GPR), which has seen much development in the statistics community, and pseudospectral approximation (PSA), which is widely used in applied mathematics and engineering. Both of these methods are routinely applied for similar purposes---replacing computationally expensive models with cheap-to-evaluate surrogates for complex analysis tasks. One of our primary goals is to compare these two approaches by analyzing the approximation spaces, experiments, and algorithms that they employ. The PSA method comes equipped with an experimental design procedure and an algorithm intended to produce approximations that are accurate in an $L^2$ sense. The GPR methodology is more flexible in that it does not come pre-equipped with an experimental design procedure. In order to compare these two methods, we will employ an experimental design criterion that also seeks accuracy in an $L^2$ sense.
%

Our main contributions are as follows. First, we develop a continuous optimization algorithm, based on sample-average approximation (SAA), to minimize an \textit{integrated posterior variance} (IVAR) design criterion for Gaussian process regression. We compare our algorithm to approaches that maximize other information-based criteria (e.g., entropy or mutual information) by evaluating their computational costs, the properties of the resulting point sets, and the accuracy of the resulting approximations.
Second, we provide a theoretical and numerical analysis comparing non-intrusive polynomial approximation---in particular pseudospectral polynomial approximation---with Gaussian process regression.  
\dedit 
While the relative performance of these methods may be a subject of broader debate, here we assess the impact of each of the three surrogate construction ingredients described above. Theoretically, we develop results to describe the difference between the approximations given the same experiments and similar approximation spaces. Numerically, we investigate the performance of the two approaches given similar approximation spaces but different experiments, chosen to be optimal for each. 

The IVAR objective\sedit, which is equivalent to the IMSE criterion~\cite{SacksWelch1989, Cressie1993}, \dedit can be minimized either over a finite (and hence discrete) design space or over a continuous design space. In the discrete case, the criterion is sometimes called the ALC (``active learning Cohn'')~\cite{Cohn1996} objective.
Minimizing the ALC criterion involves sequentially adding experiments, chosen one at a time from a discrete and finite set of candidates, to minimize a weighted average of the predictive variance. \sedit ALC has also been investigated in the context of determining \textit{local} designs for large-scale computer experiments~\cite{Gramacy2014}. \dedit ALC is often compared to other discrete design space criteria, namely ALM (``active learning MacKay'') and mutual information (MI)~\cite{Krause2008}. Minimizing the ALM criterion involves sequentially adding experiments at locations where the local predictive variance is maximized. Compared with ALM, the ALC criterion considers the effect of each experiment on the entire domain and therefore yields better designs; ALC is more expensive, however, because it requires a new variance computation for each potential design. The MI criterion sequentially seeks points that maximize the expected information gain at locations not yet chosen. MI design requires a good candidate set, which may be difficult to obtain for input domains with complex geometries, though strides have been made in this direction~\cite{Beck2014}. It also remains computationally expensive, with a complexity that grows cubically with the size of the candidate set.

Instead of dealing with the combinatorial optimization issues associated with discrete design spaces, we will pursue optimization in a continuous space. Our approach is similar to~\cite{Sacks1989,SacksWelch1989} in that we use gradient-based methods to search for optimal designs, but we will explore opportunities presented by solving the full optimization problem, without ad hoc simplifications of the design space. 
We will employ a sample-average approximation (SAA) that proves to be effective for complex domain geometries.
Benefits of this approach include generating \textit{batches} of experiments with lower computational complexity than sequentially minimizing the ALC criterion; eliminating undesirable boundary clustering effects associated with radial basis function kernels, which plague ALM designs~\cite{Ramakrishnan2005,Krause2008}; and achieving better approximation performance than either ALM or MI. 
Finally, because we perform design on a continuous space of candidate points, it becomes natural to analyze the stability and accuracy of approximation with IVAR-optimal designs from the perspective of numerical analysis. For example, in Section~\ref{sec:numerics} we will show that our algorithm generates point sets with good interpolatory properties, as measured by their Lebesgue constants. Finding these point sets via a statistical criterion raises interesting links with previous work in Bayesian numerical analysis \cite{OHagan1992,Diaconis1988}, particularly the average-case quadrature of \cite{Minka2000,Osborne2012}. A continuous design procedure also facilitates more cleanly comparing GPR and pseudospectral approximation.

Our comparison of GPR with pseudospectral approximation has two elements.
First, we use Mercer's theorem to rewrite GP regression in terms of orthogonal eigenfunctions of the kernel, such that when these eigenfunctions contain the finite basis for a pseudospectral approximation, one can directly assess the difference between the two approximations. If experimental design is based on an orthogonalizing quadrature rule, the difference between the GP mean and the pseudospectral approximation is due to eigenfunctions of the GP kernel which are not in this finite basis. These leftover eigenfunctions also account entirely for the integrated variance of the GP. 
Furthermore, we illustrate through numerical examples that experiments achieving optimal IVAR for these GPs can differ qualitatively from standard quadrature rules, and that when the IVAR criterion is then used to select experiments for GP regression, GPR can outperform pseudospectral approximation in some settings.
Second, we consider adaptive procedures for GPR that interleave IVAR-based experimental design with adaptation of the kernel hyperparameters. For test problems of moderate dimension, we find that GP approximations constructed in this way can again outperform certain adaptive pseudospectral approaches~\cite{patrick}. 

This paper is organized into two parts. The first part reviews Gaussian process regression (Section~\ref{sec:GP}), introduces the IVAR criterion and its optimization (Section~\ref{sec:expdesign}), and describes numerical comparisons of IVAR with other experimental design procedures for GP regression (Section~\ref{sec:numerics}). The second part provides a brief background on pseudospectral approximation (Section~\ref{sec:psaback}), then describes theoretical (Section~\ref{sec:theory}) and numerical (Section~\ref{sec:numCompare}) comparisons between PSA and GP regression.

%% file: gp.tex
\section{Gaussian process regression}\label{sec:GP}

Gaussian process (GP) regression \sedit can be interpreted as a Bayesian method \dedit for function approximation~\cite{OHagan1978, Rasmussen2006}, providing a posterior probability distribution over functions. The method begins with a Gaussian \sedit process \dedit prior, specified via a prior mean function $m_0(x)$ and a covariance kernel $K(x,x')$ that is positive semidefinite and bounded. Suppose that $N$ simulations of the function $f: \mathbb{R}^d \rightarrow \mathbb{R}$ are performed at parameter values $\mbf{x} := \left[x_1, \ldots x_N\right]$, $x_i \in \mathbb{R}^d$, yielding noisy function evaluations $\mbf{\hat{y}} := \left[ \hat{y}_1, \ldots \hat{y}_N\right]$, where $\hat{y}_i = f(x_i) + \xi_i$ for $i=1,\ldots,N$ and $\xi_i \sim \mathcal{N}(0, \sigma^2)$. The resulting posterior distribution is
\begin{equation*}
  \tilde{f} | \mbf{x}, \mbf{\hat{y}}  \sim\mathcal{N}(m(x), C(x,x')),
\end{equation*}
where the posterior mean is
\begin{equation}
m(x) = m_0(x) + \mbf{\alpha}^T K(\mbf{x}, x) \label{eq:postMean},
\end{equation}
and the posterior covariance is
\begin{equation}
\label{eq:postcovsimple} 
C(x,x') = K(x,x') - K(\mbf{x}, x)^T \, \mat{R} \, K(\mbf{x}, x'). 
\end{equation}
In the notation above, $K(\mbf{x}, x)$ is a (column) vector in $\mathbb{R}^{N}$ whose $i$th component is $K(x_i, x)$.  The covariance matrix $\mat{R}^{-1}$ has elements $\left[\mat{R}^{-1}\right]_{ij} = K(x_i, x_j) + \delta_{ij}\sigma^2$. Finally, the $i$th element of the vector of coefficients $\mbf{\alpha} \in \mathbb{R}^N$ is $\alpha_i = \mat{R}_{[i,:]} \left(\mbf{\hat{y}} - m_0 (\mbf{x}) \right)$. 

Gaussian process regression reverts to interpolation when $\varn=0$. However, as $N\to \infty$ the covariance matrix $\mat{R}^{-1}$ becomes ill-conditioned; a small value for $\varn$, called a nugget, is often then introduced to stabilize the procedure~\cite{Neal1997}. \sedit Note also that we have not included inference of the prior mean $m_0(x)$ in the Bayesian formulation above. If the prior mean is described via some parametric model $m_0(x) = \mbf{\beta}^T \mbf{n}(x)$, where $\mbf{n}(x)$ is a vector of basis functions, then Bayesian inference of the coefficients $\mbf{\beta}$ would add terms to the posterior covariance. In practice, however, it is common and quite effective to assume either a zero or non-zero constant term for $m_0$ and to fix its value (for example, by maximizing the log-marginal likelihood) before performing the GP update; see, e.g.,~\cite{Jeong2005,Lizotte2007,Gramacy2014}. For simplicity, we fix the prior mean here. Doing so will also help focus the comparison of nonparametric GP regression with parametric PSA in Section \ref{sec:spectral} on its essential aspects. \dedit

\subsection{Reproducing kernels and Mercer's theorem}

Many elements of this work rely on the interpretation of the covariance kernel through its eigenfunctions, and to this end we recall the properties of a Mercer kernel. Let the kernel $K:\mathcal{X} \times \mathcal{X} \to \mathbb{R}$ be defined on a first-countable \cite{Dudley2002} space $\mathcal{X} \subseteq \mathbb{R}^d$ endowed with a strictly positive Borel measure $\mu$. Suppose that the kernel is continuous, positive semi-definite,
\begin{equation}
  \quad \int_{\chi \times \chi} K(x,x') g(x) g(x') d\mu(x)d\mu(x') \geq 0, \ \forall g \in L^2_{\mu}(\mathcal{X}) \label{eq:spd},    
\end{equation}
and in $L^1(\mathcal{X}, \mu)$
\begin{equation} \label{eq:bound} 
  \int_{\chi} \left | K(x,x) \right | d\mu(x) < \infty.
\end{equation}
Additionally we define the integral operator $T_K: L^{2}_{\mu}(\mathcal{X}) \to L^{2}_{\mu}(\mathcal{X})$ such that $T_K f = \int_{\chi} K(x,x') f(x')d\mu(x').$
This operator has a countable system of eigenvalues $\lambda_j$ that are non-negative and that satisfy
\[
   \sum_{j=1}^{\infty} \lambda_j^2 < \infty.
\]
The eigenfunctions $\phi_i$ of $T_K$ form an orthonormal basis of $L^2_{\mu}(\mathcal{X})$. These eigenfunctions and eigenvalues can be used to define the reproducing kernel Hilbert space (RKHS) associated with the kernel \cite{Rasmussen2006}. Mercer's theorem lets us represent $K$ as a convergent series in terms of this eigensystem.

\begin{theorem}[Mercer]~\cite{Minh2006,Ferreira2012,Buescu2004} \label{thm:mercer} Let $\mathcal{X} \subset \mathbb{R^d}$ be first countable or locally compact, $\mu$ a strictly positive Borel measure on $\mathcal{X}$, and $K$ a continuous function on $\mathcal{X} \times \mathcal{X}$ satisfying~(\ref{eq:spd}) and~(\ref{eq:bound}). Then 
    \begin{equation}
        K(x, x') = \sum_{i=1}^{\infty} \lambda_i \phi_i(x) \phi_i(x') ,
        \label{eq:mercer}
    \end{equation}
   where the series converges absolutely for each pair $(x, x') \in \mathcal{X} \times \mathcal{X}$ and uniformly on each compact subset of $\mathcal{X}$.
\end{theorem} \\

When comparing GP regression with pseudospectral approximation in later sections of this paper, we will also use the notion of a truncated kernel. These are kernels for which $\lambda_{i > \ell} =0$, and for which we can equivalently write \eqref{eq:mercer} as
\begin{equation*}
 K(x, x') = \sum_{i=1}^{\ell} \lambda_i \phi_i(x) \phi_i(x').
\end{equation*}
One common example of a truncated kernel is a polynomial kernel, e.g., $K(x,x^{\prime}) =  (x^Tx^{\prime} + 1)^p$,  where $p$ is a positive integer.

We can use Mercer's theorem to write the integrated posterior variance of the Gaussian process in terms of the eigensystem $(\lambda_i, \phi_i)$. The posterior variance at any point in the domain is $c(x) := C(x,x)$. The integrated variance then becomes
\begin{align}
\int c(x) d \mu(x)  &= \int  \left( K(x,x) - K(\mbf{x}, x)^T \mat{R} \,  K(\mbf{x}, x)\right)  d\mu(x) \nonumber \\
          &= \int K(x,x)d\mu(x) - \int \sum_{i=1}^{\infty} \sum_{j=1}^{\infty} \lambda_i \lambda_j \pvec{i}^T \mat{R}  \pvec{j} \phi_i(x) \phi_j(x)  d\mu(x) \nonumber \\
          &= \sum_{i=1}^{\infty} \lambda_i - \sum_{i=1}^{\infty} \lambda_i^2  \pvec{i}^T \mat{R} \pvec{i}, \label{eq:intVarBasis}
\end{align}
where 
$\pvec{i} := \left[\phi_i(x_1), \ldots, \phi_i(x_N) \right]^T.$
The first term is the integrated variance of the prior. The second term, which is always non-negative, reflects reduction in the integrated prior variance due to conditioning on the data. We will now examine how experimental design procedures can use this integrated variance as an optimization objective.

%% file: expDesign2.tex
\section{Integrated variance experimental design}\label{sec:expdesign}
We will design experiments to minimize the integrated posterior variance (IVAR) of the Gaussian process. This choice is motivated by inferential considerations. As opposed to design procedures based on Latin hypercube sampling, quasi-Monte Carlo, quadrature, or other ``lattice'' designs, the present design strategy directly aims to minimize a measure of the uncertainty associated with the approximation. One advantage of this approach is that it can be used to design experiments on a wide variety of input domains, not just domains with tensor-product or some other canonical structure. Another advantage of computing and monitoring a measure of uncertainty is that it provides useful feedback about the quality of the approximation; for example, if the data do not yield much reduction in uncertainty, one can adjust the approximation space or some other aspect of the surrogate construction methodology.

\sedit
To put the IVAR criterion in context, we note that it is equivalent to an expected integrated squared error of the posterior mean.
\dedit
First consider the posterior expectation of the squared error in function values, integrated over the parameter space,
\[
\mathbb{E}_{\tilde{f}\vert \mbf{x}, \mbf{\hat{y}}} \left [\int \left( \tilde{f}(x) - f(x)\right)^2 d\mu(x) \right ],
\]
where $\tilde{f}$ can be thought of as a posterior realization of the Gaussian process and $\mu$ is the measure on the parameter space $\mathcal{X}$. 
We can divide this quantity into two terms,
\begin{align}
\mathbb{E}_{\tilde{f}\vert \mbf{x}, \mbf{\hat{y}}}  \left [\int \left( \tilde{f}(x) - f(x)\right)^2 d\mu(x) \right ] &=\int \left(m(x \vert \mbf{x}, \mbf{\hat{y}} ) - f(x)\right)^2 d\mu(x) + \int c(x \vert \mbf{x} ) d\mu(x),
\label{eq:ISEdivide}
\end{align}
where the first term on the right-hand side is the integrated squared error of the posterior mean and the second term is the integrated posterior variance, i.e., the IVAR. We have explicitly indicated all the conditioning on the right-hand side of \eqref{eq:ISEdivide}. Note that the second term is independent of the sampled values $\mbf{\hat{y}}$ of the function $f$.\footnote{In this section, we keep the prior covariance kernel $K$ fixed. In Section~\ref{sec:numCompare}, we will consider closed-loop adaptive design strategies that learn hyperparameters of the kernel $K$ from evaluations of the function $f$.}  Computing the first term, on the other hand, requires the ability to evaluate $f$. Directly using this term in a design criterion would defeat the purpose of experimental design, which is motivated by the desire to evaluate $f$ sparingly. Instead, we can consider the expectation of this squared error over the joint distribution of $f$ and $\mbf{\hat{y}}$, for a fixed design $\mbf{x}$. We assume that $f$ is drawn from the prior $\mathcal{N} \left (m_0(x), K(x, x') \right)$, and therefore this expectation becomes the Bayes risk of the posterior mean under squared error loss, which is equivalent to the integrated mean squared error (IMSE) criterion proposed by \cite{SacksWelch1989}. Some manipulation shows that this Bayes risk is indeed \textit{equal} to the integrated posterior variance, i.e., that after taking the expectation over $f$ and $\mbf{\hat{y}}$, the two terms on the right side of \eqref{eq:ISEdivide} are the same. For further details, see, e.g., \cite[p.~92]{Santner2013}.
\ssedit Thus IVAR minimization can also be understood as $I$- or $V$-optimal design \cite{Atkinson2007}, in the sense of minimizing the $\mu$-weighted variance of the predictions over the design region $\mathcal{X}$.
\ddedit

\sedit
A different connection to optimal design theory can be made by finding a  \ssedit {finite} parameterization (and hence, in general, an \textit{approximation}) of \ddedit the Mercer kernel $K$ and converting the problem into one of parametric model fitting and prediction~\cite{Fedorov1996,Fedorov1997,Fedorov2007,Harari2014}. Such approaches decompose the kernel by \ssedit computing a Karhunen-Lo\`{e}ve expansion~\cite{Fedorov2007,Harari2014} of the Gaussian process and then truncating it; \ddedit alternatively, one can use a polar spectral approximation~\cite{Spock2010} of the kernel to avoid computing its eigenfunctions. In either case, once the kernel has been decomposed and truncated, one can approximate the posterior mean of the GP prediction as a linear combination of finitely many basis functions, 
\begin{equation}
  m(x) = m_0(x) + \sum_{i=1}^{\ell} \theta_i \phi_i(x).
\end{equation}
In optimal design theory, one may then seek experiments to best learn about the mean structure $m_0$, the parameters $\theta_1, \ldots, \theta_{\ell}$, or to minimize some uncertainty in the prediction. These choices correspond to different classical ``alphabetic optimality'' criteria applied to the corresponding information matrix, e.g., $A$, $D$, $G$, \ssedit or $I$-optimal \ddedit design.

The IVAR design procedure presented in this paper, however, does not require a finite parameterization of the kernel; instead it can use a closed-form expression for the posterior covariance, with no truncation. A direct comparison between such kernel-based procedures (for example, the Sacks-Ylvisaker approach described in \cite{Fedorov1996}) and parametric optimal design theory is outside the scope of this paper; we refer readers to~\cite{Fedorov1996,Steinberg2004,Muller2005,Fedorov2007} instead. We will, however, compare our IVAR optimization procedure to other algorithms and design objectives that do not require explicitly finding a finite parameterization of the kernel. 
Later, when comparing GP regression to another surrogate modeling methodology---pseudospectral polynomial approximation, in Section~\ref{sec:spectral}---we will return to the eigenfunction viewpoint of GP regression. In that context, our focus will not be on algorithms that require explicit access to the eigenfunctions of the kernel, but rather on how the eigenfunction viewpoint exposes the distinct modeling assumptions made by the two methodologies.
\dedit

\subsection{IVAR evaluation and minimization}

For a chosen number of experiments $N \geq 1$ and a fixed prior covariance kernel $K$, our optimal experimental design is a set of evaluation points $\mbf{x}^{\ast} := \left[x_1^{\ast}, \ldots x_N^{\ast} \right]$, $x_i^{\ast} \in \mathbb{R}^d$, minimizing the integrated posterior variance:
\begin{equation}
\mbf{x}^* = \arg \min_{\mbf{x} \in \mathcal{U}} \int_\mathcal{X} c(x | \mbf{x}) d \mu(x) .
\label{eq:obj}
\end{equation}
Here $\mathcal{U}$ is the space of all feasible experiments and the posterior variance is specified via \eqref{eq:postcovsimple}.

While this objective function is similar to that in~\cite{seo2000gaussian}, here we will employ a continuous space  $\mathcal{U}$ of possible experiments. Also, we will solve the optimization problem \eqref{eq:obj} both in a non-greedy fashion (finding all $N$ design points simultaneously) and using greedy updates with varying batch sizes. In this section the number of design points $N$ and the prior covariance kernel $K$ will be considered fixed. Later, in Section~\ref{sec:numCompare}, we will consider closed-loop design procedures that alternate between batch minimization of IVAR and updates of the covariance kernel.

\subsubsection{Sample-average approximation of IVAR}\label{sec:imp1}
One method of minimizing IVAR involves numerically evaluating the objective in \eqref{eq:obj} via quadrature or a quasi-Monte Carlo or Monte Carlo (MC) sampling procedure. Since the variance is in general a smooth function of $x$, quadrature schemes may work efficiently for low-to-moderate dimensional input spaces, but Monte Carlo will generally work better in higher dimensions; Monte Carlo also offers more flexibility for non-tensor-product domains $\mathcal{X}$. 
Monte Carlo sampling replaces the integral with the summation
\begin{equation}
\label{eq:mcobjective}
\int_\mathcal{X} c(x \vert \mbf{x} )  \, d\mu(x) \approx \frac{1}{N_{mc}}\sum_{i=1}^{N_{mc}} c(\hat{x}_i  | \mbf{x} ) =: \hat{J}^{\text{mc}}(\mbf{x}),
\end{equation}
where $N_{mc}$ is the number of Monte Carlo samples and $\hat{x}_i \sim \mu$. Computing the integrated variance for a set of $N$ experiments then requires an inversion of the covariance matrix, an $\mathcal{O}(N^3)$ operation, and variance evaluations at $N_{mc}$ points, an $\mathcal{O}(N_{mc}N^2)$ operation. 

The sample-average approximation (SAA)~\cite{Shapiro1991} approach to optimization simply replaces the expectation in the objective \eqref{eq:obj} with a quadrature or Monte Carlo approximation at a \textit{fixed} set of points and minimizes this objective. After one has chosen this set of fixed points, the minimization becomes a constrained deterministic minimization problem. We can use readily available analytical derivatives of the objective in this setting. In particular, given the form of the kernel $K$, we can directly compute the gradient $\nabla_{\mbf{x}} c (x |\mbf{x})$ from \eqref{eq:postcovsimple}; details are given in Appendix~\ref{sec:gradientappendix}. The gradient of the SAA objective obtained from Monte Carlo then becomes
\begin{equation*}
\nabla_{\mbf{x}} \hat{J}^{mc}(\mbf{x}) = \frac{1}{N_{mc}} \sum_{i=1}^{N_{mc}} \nabla_{\mbf{x}} c(\hat{x}_i | \mbf{x}),
\end{equation*}
and similarly for quadrature. 
%
%
%
%

%

%

%
Finally, we note that in some cases it may be more convenient to work with closed-form expressions for the eigenfunctions of the kernel $K$ rather than the kernel itself; such situations arise when one has a desired basis of approximation but the corresponding closed-form kernel is unknown, or if the eigenfunctions can otherwise be easily computed. In this case, one can rewrite the IVAR objective in terms of eigenfunctions and simply maximize the second term on the right of~\eqref{eq:intVarBasis}. Written in this form, the objective does not require integration with respect to the parameter measure $\mu$. Indeed, having the eigenfunctions in hand is tantamount to already having performed the integration, as the eigenfunctions are solutions of the homogeneous Fredholm integral equation with operator $T_K$. \cite{Gauthier2014} uses this approach for IVAR minimization, exploring truncations of \eqref{eq:intVarBasis} to $\ell$ eigenfunctions.

\subsubsection{Batch and greedy implementations}\label{sec:imp3}
\ssedit In this work, we have used gradient-based optimization algorithms from both the NLopt~\cite{nlopt} and SciPy~\cite{Jones2014} optimization packages, with similar performance,
 to solve the optimization problem \eqref{eq:obj}. \ddedit This problem involves $N$ design points, each in $d$ dimensions, and thus has $Nd$ unknowns. 
It can become expensive to solve when $Nd$ is very large. In these cases it may be useful to solve a sequence of smaller optimization problems to achieve an approximate solution. In the numerical examples below, we will investigate constructing these smaller problems through the use of a greedy minimization procedure. In this procedure one decides how many training points $M$ are computationally feasible to minimize. Suppose that $k=N/M$ is an integer. Then the greedy procedure solves $k$ optimization problems of size $Md$. Once the $(j-1)$th problem is solved, for $j \leq k$, we have obtained the experiments $\mbf{x}_{j-1}$. During the $j$th iteration we find $M$ points to append to $\mbf{x}_{j-1}$.

\subsection{Entropy and mutual information}
Statistical criteria underlie many other experimental design procedures for Gaussian process regression \cite{Santner2003}. Two popular techniques include minimizing the conditional entropy of the Gaussian process at unobserved locations, or maximizing the mutual information (MI) between the locations at which experiments are performed and the rest of the design space. Our numerical results will compare IVAR designs with MI and entropy-based designs because the latter have algorithms specifically tailored for Gaussian process regression. For a comparison between these two methods and additional design procedures, e.g., based on classical alphabetic optimality criteria, we refer to~\cite{Krause2008}.

The conditional entropy design procedure seeks experiments that reduce the uncertainty, as measured by entropy $H$, across a typically finite set of possible simulation locations. If the set of candidate experimental locations is denoted by $\mathcal{D}$ and $\mbf{x} \subset \mathcal{D}$ are the chosen locations, then $x^c := \mathcal{D} \setminus \mbf{x}$ are the locations at which the entropy is evaluated. In particular, one seeks $\mbf{x}$ to minimize the conditional entropy \sedit $H( F_{\mbf{x}^c} | F_{\mbf{x}} ) = -\int p(f_{\mbf{x}^c},f_{\mbf{x}}) \log p(f_{\mbf{x}^c}|f_{\mbf{x}}) d f_{\mbf{x}^c} d f_{\mbf{x}}$, where $F_{\mbf{s}}$ is a random variable representing the outputs of the simulation model at a set of inputs $\mbf{s}$ and $p$ denotes a probability density. \dedit Minimizing this function is shown to be NP-hard in~\cite{Ko1995}. In practice, one instead employs greedy but suboptimal algorithms that add one experiment at a time---for example, adding each experiment at the location where the conditional entropy is largest~\cite{Shewry1987,Cressie1993, MacKay1992}. In the context of GPR,  $F_{\mbf{x}}$ and $F_{\mbf{x}^c}$ are jointly Gaussian; this procedure then becomes equivalent to greedily adding experiments at the locations of highest variance, and is commonly called the MacKay criterion (ALM). Other algorithms for choosing a subset of points to minimize the conditional entropy are based on the DETMAX algorithm~\cite{Mitchell1974}, demonstrated for GP regression in~\cite{Currin1991} 

The mutual information criterion for experimental design considers the change in the entropy of $F_{\mbf{x}^c}$ before and after performing experiments at locations $\mbf{x}$: $H (F_{\mbf{x}^c} ) - H( F_{\mbf{x}^c} \vert F_{\mbf{x}} )$, which is equivalent to the mutual information of $F_{\mbf{x}}$ and $F_{\mbf{x}^c}$. This objective is also typically maximized in a greedy fashion: from the candidate set of experiments $\mathcal{D}$, the element which yields the greatest MI is chosen at each iteration. Unlike the conditional entropy, however, MI is submodular, guaranteeing that a greedy approach performs within a constant factor of the full $N$-experiment maximization \cite{Krause2008}. The greedy procedure requires the inversion of a matrix containing the covariance between every pair of candidate experiments, which arises when computing the entropy on the set of all simulation locations. If one is choosing $N$ experiments from a set of $M \gg N$ candidates, each iteration then requires inverting a matrix of size $M \times M$. This expense is typically large and many recent efforts have aimed at reducing it. \cite{Krause2008} reduces the cost of each iteration to $\mathcal{O}(M)$ by using specialized local kernels. \cite{Beck2014} points out that the quality of an MI design depends crucially on the candidate set, and modifies the greedy algorithm of \cite{Krause2008} to be more robust by resampling the set of candidate experiments after each new point is chosen.

Besides the choice of objective, our IVAR-based design algorithm differs from the entropy and MI approaches above in other fundamental ways. First, we select experiments from a continuous design space and thus avoid the challenges of combinatorial optimization. In this sense, we follow the approach of~\cite{SacksWelch1989} by solving a continuous optimization problem using gradient descent algorithms. Second, as described in Section~\ref{sec:imp3}, our approach can identify multiple experiments simultaneously. Designing batch experiments is advantageous because interactions between the experiments are taken into account; we will demonstrate this advantage empirically in the next section. 
Our continuous design approach also takes several steps beyond~\cite{SacksWelch1989}. First, we do not employ particular patterns or partitions to simplify the design space. And as described earlier, we introduce a sample-average approximation of the IVAR objective. This helps design experiments on complex input domain geometries by penalizing proximity to the domain boundaries; the objective automatically becomes large in locations where there are few samples. We will also address the ``pileup'' problem identified in~\cite{SacksWelch1989}, explaining it in terms of the design size relative to correlation kernel complexity.

\subsection{Computational complexity}
In some sense, if evaluating $f$ is sufficiently expensive, then the computational cost of finding a good design is immaterial. But the design procedures described above can have very different costs, and in practice one may not want the computational effort required for experimental design to be too large. Table~\ref{tab:compComplex} summarizes the computational complexity of evaluating and optimizing the various experimental design objectives considered above. In the IVAR scenarios, $N_{mc}$ denotes the number of Monte Carlo points sufficient for \eqref{eq:mcobjective}, where typically $N_{mc} \gg N$. For the entropy and MI criteria we assume that the number of candidate designs $N_{E} = | \mathcal{D} | \gg N$. Finally, we assume that optimizing the IVAR objective requires $L$ objective and/or gradient evaluations, where typically $L \ll N_{mc}$.

Greedy minimization of the conditional entropy (ALM) requires $N$ iterations. In each iteration the variance must be computed at $N_{E}$ locations.\footnote{Although we have described a discrete approach for ALM minimization, a continuous optimization approach is also possible. But the pointwise variance of a GP has many local maxima, and a multi-start procedure would likely be necessary for continuous optimization to be competitive with discrete enumeration. It is then unclear whether continuous optimization would in fact be more efficient in this case.} At iteration $k< N$, this variance computation requires an $\mathcal{O}(k^3)$ inversion followed by $N_E$ variance evaluations, each of complexity $\mathcal{O}\left(k^2\right)$. Thus the complexity for step $k$ is $\mathcal{O}(k^3 + N_E k^2)$. The total complexity is thus $\mathcal{O}\left( \sum_{k=1}^N k^3 + N_E k^2 \right) = \mathcal{O} \left( N^4 + N_E N^3 \right)$.

\begin{table}[h!]
\center
\begin{tabular}{|c|c|c|c|}
\hline
Design objective & Optimization method & Objective eval & Optimization cost \\[2pt]
\hline
\multirow{2}{*}{IVAR} & Monte Carlo, batch  & $\mathcal{O}\left(N_{mc}N^2 + N^3\right)$ & $\mathcal{O}\left( L N_{mc}N^2 + L N^3 \right) $\\
& eigenfunction form, batch & $\mathcal{O}\left( \ell N^2 + N^3 \right)$ & $\mathcal{O}\left( \ell L N^2 + L N^3 \right)$ \\[2pt]
\hline
entropy &  greedy & $\mathcal{O}(N^3 + N_E N^2)$ & $\mathcal{O}(N^4 + N_{E}N^3)$ \\[2pt]
\hline
\multirow{1}{*}{mutual information} & greedy & $\mathcal{O}(N_{E}^4)$ & $\mathcal{O}(N N_{E}^4)$ \\[2pt]
\hline
\end{tabular}
\caption{Computational complexity of different experimental design algorithms.}
\label{tab:compComplex}
\end{table}

The ALM approach is often the computationally cheapest option, but not if $N_E > L N_{mc}/N$. If $N_E$ and $N_{mc}$ are of comparable magnitude, then the comparison of entropy and IVAR minimization depends on whether $L/N < 1$, i.e., how the number of optimization iterations (in the continuous case) compares to the number of candidate design points (in the discrete case). 
We also see that the MI design becomes very expensive for large $N_{E}$; a small candidate design set must be chosen for the MI procedure to remain tractable. Using a smaller candidate set, however, increases the possibility that the chosen designs will perform poorly.

%% file: numExamples.tex
\section{Numerical examples of IVAR designs}~\label{sec:numerics}
We now provide numerical examples illustrating the quality of designs arising from continuous IVAR minimization. 
All of the examples in this section were performed using the freely available GPL-licensed GPEXP package~\cite{gpexp} for python. 
First, we examine the interpolatory properties of IVAR design points. Then we compare IVAR, entropy, and MI-based designs on several domains.
\ssedit
We also include comparisons with space-filling designs obtained by the method of Lekivetz and Jones (LJ)~\cite{Lekivetz2015}, which is well suited for many non-rectangular regions. The LJ algorithm works by grouping randomly sampled points into equally sized clusters using Ward's mininum variance criterion; one design point is then extracted from each cluster (for instance, by computing the cluster centroid). The LJ approach encounters difficulties for non-convex domains, and thus we do not use it in the non-convex test case of Section~\ref{sec:mickey}.
\ddedit

\input{rbf_interpolation2.tex}

\subsection{Approximation comparisons} \label{sec:irregularDomains}
Here we illustrate IVAR designs on variety of input domains $\mathcal{X}$, and compare the performance of IVAR designs with that of designs obtained through conditional entropy minimization or MI maximization. \sedit We highlight designs and approximations on irregular input domains (e.g., domains that are neither hypercubes nor $\mathbb{R}^d$), which are critical in many real-world applications. In particular, irregular domains often arise as a result of domain partitioning by a discontinuity detection algorithm, for models whose outputs are piecewise smooth. For example, in~\cite{Jakeman2011,Gorodetsky2014} the authors automatically partition the input domain of a genetic toggle switch model that exhibits a phase transition. Following this partitioning, function approximation proceeds on two separate but irregular domains. In~\cite{Sargsyan2009,Sargsyan2012} the authors study a climate model which exhibits a discontinuity. Discontinuity detection is used to split the input domain into two irregular subdomains, and is followed by function approximation. We will thus attempt to show the applicability of our algorithm to input domains that are characteristic of such problems, which have no guarantees of convexity or even connectedness. \dedit

In most of our numerical experiments, we employ an isotropic squared exponential kernel $K(x,y) = \exp(-\Vert x-y \Vert^2/2l^2)$ and we set $\varn=10^{-10}$; \ssedit the only exception is in Section~\ref{sec:periodic}, which uses a periodic kernel to contrast space-filling designs with designs that are adapted to the kernel. \ddedit Entropy minimization is pursued through the ALM approach described earlier; each experiment is chosen by comparing the variance at $N_E=10^4$ possible experimental locations sampled according to the parameter measure $\mu$ on $\mathcal{X}$. \ssedit For MI maximization, we select designs from $N_E=200$ candidate locations generated by the space-filling LJ design. \ddedit The size of the MI candidate set is chosen so that the computational times of the different experimental design procedures are comparable; it is also comparable to sizes used in the literature~\cite{Beck2014}. We find that small changes in the size of this candidate set do not qualitatively change the results reported below. 
\ddedit

Finally, we define the relative $L^2$ error of a GP approximation as $\Vert f - m \Vert _{L^2_{\mu}} / \Vert f \Vert _{L^2_{\mu}}$, and we estimate it using $10^5$ Monte Carlo samples. 
\sedit 
While in the previous section we considered $L^{\infty}(\mathcal{X})$ error, recall that the IVAR objective function reflects the expectation of a squared $\mu$-averaged error. Thus $L^2_{\mu}$ is a logical metric of quality. Other efforts, e.g.,~\cite{Loeppky2010}, empirically evaluate the $L^2$ and $L^{\infty}$ errors of approximations resulting from a variety of other design criteria.
\dedit
\subsubsection{Circular domain}\label{sec:circle}
\begin{figure}
\centering
\includegraphics[scale=1.0]{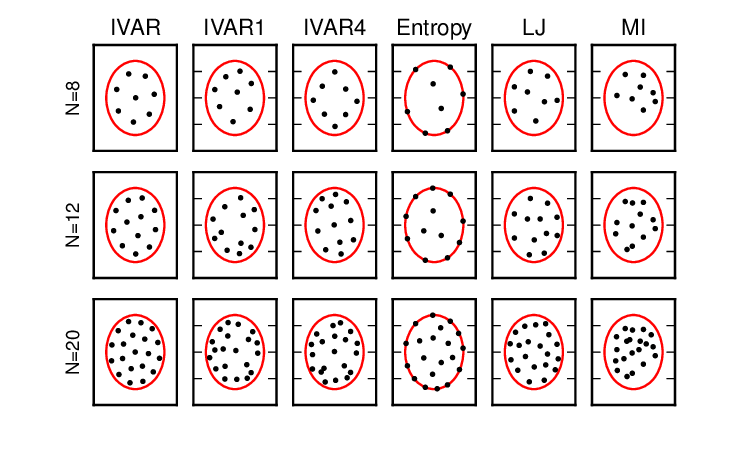}
\caption{Experimental designs on a circuar domain obtained using various strategies for an isotropic squared exponential kernel with $l=0.2$. IVAR-$M$ corresponds to a greedy IVAR strategy, adding points in batches of $M$.}
\label{fig:circle20}
\end{figure}

We first construct designs on a circular domain in $\mathbb{R}^2$ with radius 0.7. We fix the correlation length in the kernel $K$ to $l=0.2$; the measure $\mu$ is uniform over the domain. Results from each design strategy are shown in Figure~\ref{fig:circle20}, for designs with $N=8$, $12$, and $20$ points. We consider full batch IVAR, designing all $N$ points simultaneously, along with greedy IVAR strategies that add points in groups of $M=1$ or $M=4$. Full IVAR results in the most symmetric and regularly spaced points. But all of the IVAR strategies attempt to spread the points throughout the domain, whereas entropy design first places points on the boundaries, which is not a desirable feature for radial basis function kernels~\cite{Ramakrishnan2005,Krause2008}. \ssedit The LJ designs are reasonably space-filling, as expected, but full IVAR yields even more consistent spacing.
\ddedit
 For timing reference, we note that finding a 40-point design with the full IVAR criterion took 11 seconds, IVAR-1 took 65 seconds, IVAR-4 took 29 seconds, entropy-based design took 10 seconds, and MI design took 79 seconds. Thus non-greedy IVAR design and entropy design take approximately the same amount of time, while MI maximization is slightly more costly, though of the same order of magnitude. 

\begin{figure}
\centering
\subfigure[Two-dimensional domain; prior correlation length $l=0.2$]{
	
	\includegraphics[scale=0.8]{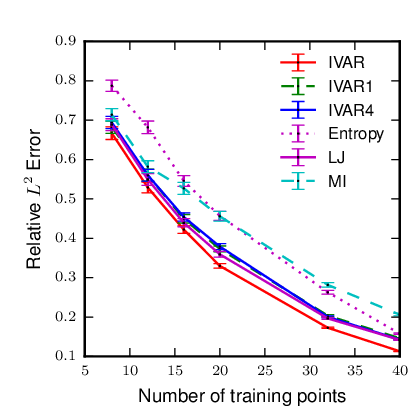}
	\label{fig:circleNodesApprox}
}
\subfigure[Five-dimensional domain; prior correlation length $l=0.5$] {
	\includegraphics[scale=0.8]{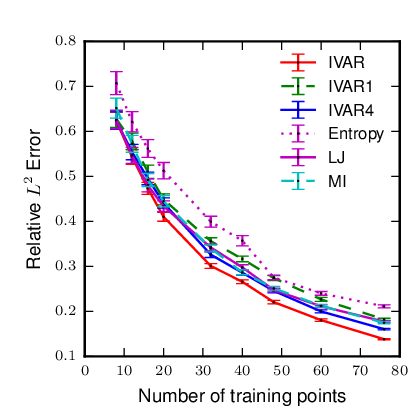}
    \label{fig:5dcircleNodesApprox}
}
\caption{Relative $L^2$ approximation error on 1000 functions drawn from a GP with squared exponential kernel; error bars indicate sample standard deviations of these errors. Domains $\mathcal{X}$ are balls in $\mathbb{R}^2$ or $\mathbb{R}^5$.}
\label{fig:spherical experiments}
\end{figure}

To evaluate the effectiveness of these designs for function approximation, we perform GP regression on 1000 functions independently sampled from the prior Gaussian process, with results shown in Figure~\ref{fig:circleNodesApprox}. For each sampled function $f$, we compute the relative $L^2$ error between the posterior mean $m$ of the GP and $f$, as described above. We then report the average and standard deviation of this error, for each design strategy and different values of $N$. The IVAR designs, including the IVAR-$M$ greedy strategies, clearly outperform the entropy and MI designs. Optimality of the full IVAR designs is to be expected, since we are essentially calculating the Bayes risk (the expectation of $\Vert f - m \Vert^2_{L^2_\mu}$) which is equivalent to IVAR, as discussed in Section~\ref{sec:expdesign}. But it is noteworthy that even the greedy IVAR strategies show better performance than the MI and entropy designs. 
\ssedit
The LJ designs show errors comparable to the greedy IVAR strategies, but not as low as full IVAR.
\ddedit
In Figure~\ref{fig:5dcircleNodesApprox}, we repeat this study for a higher-dimensional domain---a ball in $\mathbb{R}^5$---and observe similar trends.
\sedit
For reference, the computational times required to find 40-point designs in $d=5$ dimensions are: 13 seconds for batch IVAR, 125 seconds for IVAR-1, 45 seconds for IVAR-4, 10 seconds for entropy, and 82 seconds for MI. These results further support the idea that full IVAR is computationally competitive with entropy and that design with MI is more expensive.
\dedit

\subsubsection{Periodic kernel}\label{sec:periodic}
\begin{figure}
\centering
\includegraphics[scale=1.0]{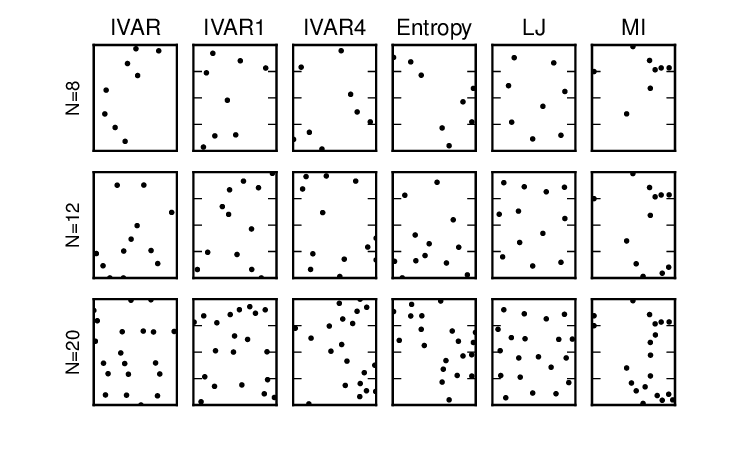}
\caption{\ssedit Experimental designs obtained using various strategies for a periodic kernel on the domain $[-1,1]^2$.\ddedit}
\label{fig:pernodes}
\end{figure}

\begin{figure}
\centering
\includegraphics[scale=0.8]{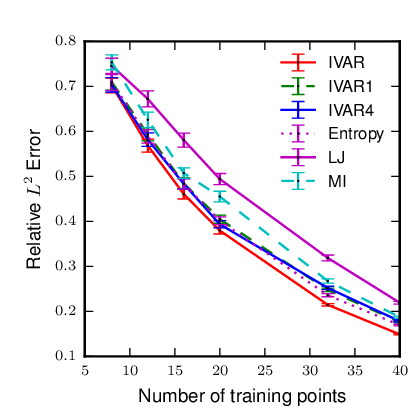}
\caption{\ssedit Relative $L^2$ approximation error on 1000 functions drawn from a GP with periodic kernel on $[-1,1]^2$ with parameters $p=1.0$ and $l = 0.9$; 
error bars indicate sample standard deviations of these errors.\ddedit}
\label{fig:period}
\end{figure}

\ssedit
The IVAR, MI, and entropy-based design criteria explicitly incorporate the kernel of the GP, while the space-filling LJ design approach does not. Thus it is reasonable to expect kernel-adapted designs to be more successful for a broader range of kernel specifications. To assess this behavior, we repeat the experiments of the previous subsection on a square domain with a periodic kernel. In particular, we use the kernel
\begin{equation*}
K(x,y;p,l) = \exp\left( -2 \sum_{i=1}^d \frac{\sin^2( \pi |x-y| / p)}{l^2}\right),
\end{equation*}
on the domain $[-1,1]^2$, with period $p=1.0$ and correlation length $l=0.9$. The resulting design points are shown in Figure~\ref{fig:pernodes} and the function approximation results are shown in Figure~\ref{fig:period}. Now, optimal designs for the IVAR, entropy, and MI schemes are {not} space filling. They are clustered in areas that reflect the periodic structure of the covariance kernel. The space-filling LJ nodes do not exploit this property and, as shown in Figure~\ref{fig:period}, yield higher errors than the kernel-adapted nodes.

\ddedit

\subsubsection{Non-convex, non-simply connected domain}\label{sec:mickey}
\begin{figure}
\centering
\subfigure[12 points]{
	\includegraphics[scale=0.8]{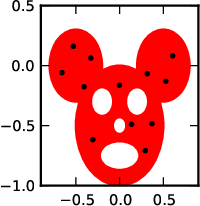}
    \label{fig:mickey12pts}
} 
\subfigure[16 points]{
	\includegraphics[scale=0.8]{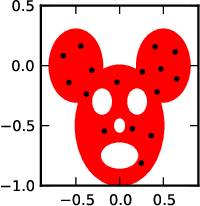}
    \label{fig:mickey16pts}
} 
\subfigure[32 points]{
	\includegraphics[scale=0.8]{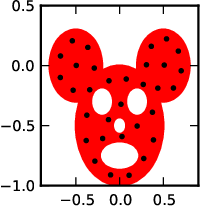}
    \label{fig:mickey32pts}
}
\subfigure[60 points]{
	\includegraphics[scale=0.8]{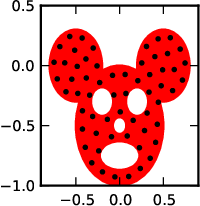}
    \label{fig:mickey60pts}
}
\caption{Experimental designs computed by minimizing the IVAR cost function over a non-convex and non-simply connected domain; the domain is the red/shaded region, endowed with a uniform parameter measure.}
\label{fig:mickeyNodes}
\end{figure}

\begin{figure}
\centering
	\includegraphics[scale=1.0]{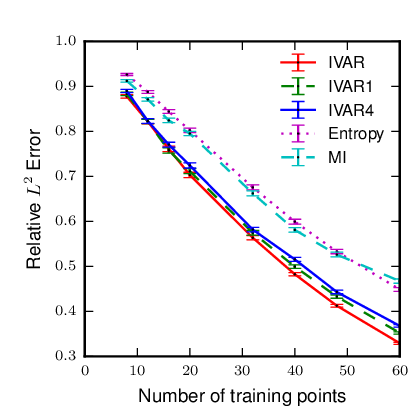}
	\label{fig:2dmickeyapprox}
\caption{Performance of IVAR minimization on non-convex, non-simply connected domain with $l = 0.1$.}
\label{fig:mickeyCost}
\end{figure}

Figure~\ref{fig:mickeyNodes} illustrates experimental design on a parameter domain that is neither convex nor simply connected, found via full (non-greedy) IVAR minimization. These designs are obtained with an isotropic squared exponential kernel with correlation length $l=0.1$. The IVAR objective is minimized using the SAA approach described in Section~\ref{sec:imp1}. Optimization for each $N$-point design began from a single randomly-generated point in the $2N$-dimensional design space; no multi-start procedure was used here.
The designs do a good job covering the domain. While twelve- and sixteen-point designs are not completely evenly distributed, designs using 32 and 60 points are very well spaced and have interesting symmetries. The approximation performance of IVAR designs---even greedy IVAR-$M$ designs---is also superior to that of the entropy and MI approaches, as shown in Figure~\ref{fig:mickeyCost}. 

%% file: rbf_interpolation2.tex
\subsection{Stability of interpolation}

As discussed above, GP regression with $\varn=0$ (and a kernel rank $\ell > N$) corresponds to interpolation. In this setting, it is useful to analyze the quality of interpolation on IVAR design points. The Lebesgue constant is often used to bound the error of a polynomial interpolant relative to the best approximation in a polynomial space of equivalent degree, but it can also be used to analyze interpolation with positive definite kernels, corresponding here to the posterior mean of the GP. Let $m(x)$ denote the interpolant or posterior mean, and let $m^\ast(x)$ denote the \textit{best} approximation of $f$, in the $L_\infty$ sense, from the finite-dimensional kernel space $H(K, \mbf{x}_N) = \Span \{ K(\cdot, x_1), \ldots, K(\cdot, x_N) \}$. We restrict our attention to approximation on compact domains $\mathcal{X} \subseteq \mathbb{R}^d$. Then the $L^\infty$ interpolation error is bounded as \cite{Fasshauer2011}
\begin{equation*}
\Vert f - m \Vert_{L_\infty(\mathcal{X})} \leq (1 + \Lambda_{K,\mbf{x}_N}) \Vert f - m^\ast \Vert_{L_\infty(\mathcal{X})},
\end{equation*}
and moreover we have \cite{de2010stability},
\begin{equation*}
\Vert m \Vert_{L_\infty(\mathcal{X})} \leq \Lambda_{K,\mbf{x}_N} \Vert  \mbf{y} \Vert_{\ell_\infty} .
\end{equation*}
The Lebesgue constant $\Lambda_{K,\mbf{x}_N}$ for prior kernel $K$ and design points $\mbf{x}_N$ is
\begin{equation}
 \Lambda_{K,\mbf{x}_N} = \max_{x \in \mathcal{X}} \sum_{j=1}^N |u_j(x)|,
\end{equation}
where $\{ u_j(x) \}_{j=1}^N$ are the \textit{cardinal functions} satisfying $u_j(x_i) = \delta_{ij}$, such that the interpolant can be written as $m(x) = \sum_{j=1}^N y_j u_j(x)$ with $y_j = f(x_j)$. We evaluate the cardinal functions as in~\cite{Fasshauer2011} by solving the linear system
\begin{equation}
\mat{R}^{-1} \mbf{u}(x) = K(\mbf{x},x),
\end{equation}
where $\mbf{u}(x) = [u_1(x) \ldots  u_N(x)]$. 

\begin{figure}
\begin{center}
\includegraphics[scale=0.8]{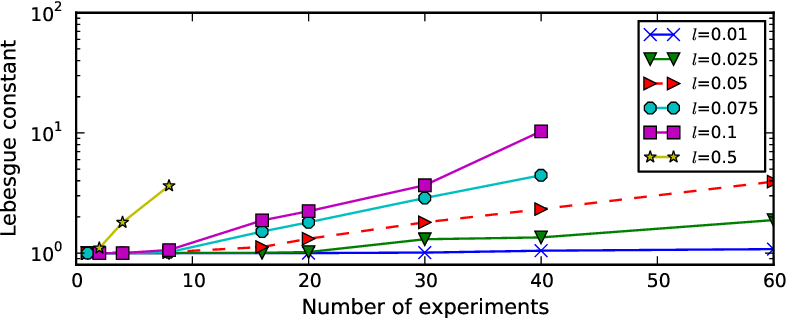}
\caption{Lebesgue constants for $N$-point IVAR designs on $[-1,1]$ with a squared exponential kernel and various correlation lengths $l$.}
\label{fig:LebesgueConstants}
\end{center}
\end{figure}

\begin{figure}
\begin{center}
\includegraphics[scale=0.8]{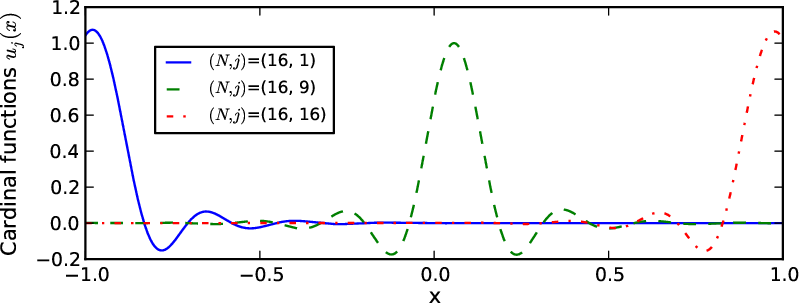}
\caption{Cardinal functions for interpolation on 16 IVAR points with a squared exponential kernel and $l=0.1$}
\label{fig:cardinal}
\end{center}
\end{figure}

Now we conduct a simple numerical experiment to evaluate the Lebesgue constants of point sets arising from IVAR minimization. We find IVAR designs on the domain $\mathcal{X} = [-1,1]$ with squared exponential kernel $K(x,x^{\prime}) = \exp({- (x-x^{\prime})^2/ 2 l^2})$. Figure~\ref{fig:LebesgueConstants} shows the associated Lebesgue constants as a function of number of design points $N$, for various correlation lengths $l$. 
Several interesting trends are observed. First, we see that for any value of $l$, the Lebesgue constant is exactly one for sufficiently small $N$. This observation is consistent with the asymptotic estimate for $l \rightarrow 0$ in \cite{RiemenschneiderS1999}. The Lebesgue constant reverts to one for small point sets because the IVAR criterion ensures that the points are well separated, and hence the cardinal functions do not overlap. Once a certain threshold value of $N$ is attained, however, the Lebesgue constant begins to increase; this threshold value is smaller for larger values of the correlation length $l$. This transition coincides with interactions among the cardinal functions: Figure~\ref{fig:cardinal} shows cardinal functions for $N=16$ and $l=0.1$, just beyond the threshold where interaction becomes significant. In this regime, the Lebesgue constant increases steadily with $N$. For a sufficiently large $N$, $\mat{R}^{-1}$ becomes poorly conditioned and direct computations of the interpolant, the cardinal functions, and the Lebesgue constant are no longer numerically stable. From the Bayesian perspective, this ill-conditioning corresponds to the ``complexity'' of the RKHS associated with the prior kernel---i.e., the effective number of nonzero eigenvalues in \eqref{eq:mercer}---being exceeded by the data.

The success of the IVAR optimization procedure itself also depends on how $N$ relates to the complexity of the kernel. In the small-$N$ and small-$l$ regime (intuitively, ``too few'' points relative to the complexity of the kernel), the IVAR cost function is relatively flat as a function of the design coordinates $\mbf{x}$, and distinguishing the quality of different designs becomes difficult. By contrast, in the large-$N$ and large-$l$ regime (``too many'' points relative to the complexity of the kernel), the IVAR value itself is exceedingly small and the problem is poorly conditioned, as described above. Away from these limiting regimes, however, the IVAR design procedure yields relatively slow growth in the Lebesgue constant and thus relatively stable interpolation.

%% file: psa.tex
\section{Comparisons with pseudospectral approximation}\label{sec:spectral}

\sedit
Pseudospectral approximation (PSA) is a well-established approach for computing polynomial approximations from pointwise function evaluations $f(x)$ \cite{Xiu2009,Constantine2012,patrick}. Recall that PSA comes equipped with an experimental design procedure and an algorithm which seeks accuracy in an $L^2_{\mu}$ sense. In Sections \ref{sec:GP}--\ref{sec:numerics}, we described an experimental design procedure for GP regression that also targets an expected $L_{\mu}^2$-type error. Thus we are now well positioned to make a comparison between GP regression and PSA. In particular, our comparison of these surrogate construction methods arises from the perspective of the three choices described in the introduction: the approximation space, the experiments, and the algorithm for constructing the approximation itself. Our comparison comprises three components described in Sections~\ref{sec:theory},~\ref{sec:spectra}, and~\ref{sec:numCompare}. We first discuss a theoretical result bounding the error between the GP posterior mean and the PSA under orthogonalizing designs. Then we show what happens when we relax the idea of exact orthogonality. Finally, we provide some numerical comparisons between the two methods on canonical approximation problems.

Overall, our goal is not so much to compare performance but rather to understand the links and distinctions between these approaches. In particular, we would like to know: are GP regression and pseudospectral approximation ever equivalent? Do the point sets used for pseudospectral approximation yield low IVAR when used for GP regression? How do IVAR-minimizing designs compare with quadrature rules, for appropriate choices of kernel?
\dedit

\subsection{Pseudospectral approximation}~\label{sec:psaback}
\sedit In this section, we recall the basics of pseudospectral approximation. \dedit Consider a set of basis functions $\{\psi_i(x)\}$ comprising a complete orthonormal system in $L^2_{\mu}$. Pseudospectral approximation takes advantage of orthonormality
\begin{equation} \label{eq:ortho}
\left< \psi_i, \psi_j \right>_{\mu} = \int \psi_i(x) \psi_j(x) d\mu(x) = \delta_{ij}, \nonumber
\end{equation}
to compute the projection of a function $f$ onto a subspace of $L^2_\mu$ spanned by $\ell$ basis functions.  The orthogonal projection
\begin{equation}
f_{\ell} (x)= \sum_{i=1}^{\ell} \left< f, \psi_i\right>_{\mu} \psi_i(x). \nonumber
\end{equation}
converges in the $L^2_{\mu}$ sense as the subspace grows: $\lim_{\ell \to \infty} \int \left(f_{\ell} - f \right)^2 d\mu = 0$.
Pseudospectral approximation departs from exact orthogonal projection by using numerical quadrature to approximate the inner products $\left<f,\psi_i\right>$. In particular, one seeks an \textit{orthogonalizing} set of nodes and weights, $\mathcal{Q} = \left\{(x_k, w_k) : k = 1 \ldots N\right\}$, i.e., a quadrature rule that computes the inner products between the first $\ell$ basis functions exactly:
\begin{equation}
\left< \psi_i, \psi_j \right>_{\mu}  = \sum_{k=1}^N w_k \psi_i(x_k) \psi_j(x_k), \ i, j=1,\ldots, \ell.
\label{eq:orthoQuad}
\end{equation}
This choice eliminates internal aliasing error \cite{patrick}, though the pseudospectral approximation will still exhibit external aliasing error, since $f$ is not necessarily in the span of $\{ \psi_1, \ldots, \psi_\ell \}$.
Given an orthogonalizing rule \eqref{eq:orthoQuad}, a pseudospectral approximation $\hat{f}_{\ell}$ can be written as
\begin{equation}\label{eq:regularpsa}
\hat{f}_{\ell}(x)  = \sum_{i=1}^{\ell} \left( \sum_{k=1}^N w_k f(x_k) \psi_i(x_k)\right)\psi_i(x) = 
\sum_{k=1}^N \left[  w_k f(x_k) \left( \sum_{i=1}^{\ell}   \psi_i(x_k) \psi_i(x) \right) \right].
\end{equation}
In subsequent analysis, it will be convenient to express $\hat{f}_{\ell}(x)$ in matrix notation. 
Let $\mathbf{W} := \diag(w_1, \ldots w_N)$, $\mbf{y} := \left[ f(x_1), \ldots, f(x_N)\right]$  and
$\psvec{i} := \left[\psi_i(x_1), \ldots, \psi_i(x_N) \right]$. Then 
\begin{equation}\label{eq:matpsa}
\hat{f}_{\ell}(x) =  \mbf{y}^T \mathbf{W} \sum_{i=1}^{\ell}\psvec{i}\psi_i(x) .
\end{equation}

Examples of $\{\psi_i\}$ and $\mathcal{Q}$ include one-dimensional orthogonal polynomial families and their corresponding Gaussian quadratures; or tensorized versions of each in multiple dimensions. However, the basis functions $\psi_i$ do not in general need to be polynomials, and other quadrature rules besides Gaussian rules may exist. Using Mercer's theorem, we can already see that $\sum_{i=1}^{\ell}\psi_i(x)\psi_i(y)$ can be interpreted as a truncated kernel $K^\ell(x,y)$ with eigenvalues that do not decay, and we can interpret the pseudospectral approximation given in~\eqref{eq:matpsa} as a weighted combination of such kernels.

%% file: theory.tex
\subsection{Same approximation space and experiments}\label{sec:theory}

To begin our comparison, we will fix two of the choices involved in surrogate modeling and evaluate the impact of the third. In particular, this section will compare the impact of different \textit{algorithms}---i.e., using pseudospectral approximation versus GP regression---for identical experiments and for the same or similar approximation spaces. Because pseudospectral approximation requires experiments to be chosen according to an orthogonalizing rule, we will also use these experiments for GP regression. We are now left to relate the basis for the pseudospectral approximation to the GP kernel, as these objects determine the approximation space.

Theorem~\ref{thm:mercer} lets us represent any Mercer kernel, and the corresponding GP posterior mean function, using the eigenfunctions.  With this connection, we first develop a more general result than required. Specifically, we show that when the basis for pseudospectral approximation comprises a subset of the eigenfunctions of a given kernel, then the $L^2_{\mu}$-norm of the difference between the resulting GP posterior mean function and the pseudospectral approximation may be bounded. Results for identical approximation spaces follow immediately as a corollary of this general case.

For the result below, we will assume that we have a Mercer kernel consisting of $\ell_{GP}$ eigenfunctions, where $\ell_{GP}$ could be finite or infinite, and a spectral expansion consisting of the first $\ell \leq \ell_{GP}$ eigenfunctions, where $\ell$ is finite. The approximations will be constructed from $N$ evaluations of the function $f$, performed at nodes of a quadrature rule $\left\{ (x_i, w_i) \right\}, i = 1 \ldots N$, that orthogonalizes the first $\ell$ eigenfunctions. Clearly $N$ depends on $\ell$. 
Let $\mat{USU}^T$ be the eigendecomposition of the matrix of covariances $K(x_i, x_j)$ among all the design points computed without a nugget, i.e., $\mat{R}^{-1} - \sigma^2\mat{I} = \mat{USU}^T$. \sedit $\mat{S}$ is a diagonal matrix with elements $\mat{S} = \diag(s_1,\ldots,s_N)$. \dedit Finally, assume a zero mean prior $m_0(x) = 0$ for the GP model. 
\sedit
The latter assumption does not restrict the generality of our results. In fact, one can transform the problem from one of approximating $f$ to one of approximating $g = f-m_0$, for some fixed function $m_0(x)$. Of course, this transformation requires one to translate the data $\mbf{y}$ accordingly. 
\dedit
Under these conditions, we have the following result, whose proof is given in Appendix~\ref{sec:errApp}. 

\smallskip

\begin{theorem}\label{thm:GPS}
    Let $\hat{f}_{\ell}(x)$ be a pseudospectral approximation represented with basis functions $\left\{\psi_1, \ldots, \psi_{\ell}\right\}$, computed via an orthogonalizing quadrature rule $\mathcal{Q}_{\ell}$ as in \eqref{eq:orthoQuad}--\eqref{eq:regularpsa}. Let \sedit $M = \underset{(x,w) \in \mathcal{Q}_{\ell}, \ i \in \{ 1, \ldots, \ell\}}{\max} \vert \psi_i(x) \vert$.  Also, let $w_{\max} = \underset{i \in \{1, \ldots, N\}}{\max} w_i$ \dedit. Let $m(x)$ be the posterior mean of GP regression with prior covariance kernel $K(x,x^{\prime}) = \sum_{i=1}^{\ell_{GP}} \lambda_i \phi_i(x)  \phi_i(x^{\prime})$ and nugget term $\sigma^2 > 0$, constructed from function evaluations at the nodes of $\mathcal{Q}_\ell$. If $\ell, N < \infty$, $\psi_i = \phi_i$ for $i=1 \ldots \ell$, and $\ell \leq \ell_{GP}$, then the difference between these two approximations is bounded as
\begin{equation}
    \sedit
    \Vert m - \hat{f}_{\ell}\Vert_{L^2_{\mu}}^2 \leq \Vert \mbf{y} \Vert_2^2 \left(\frac{\ell N M^2 w_{\max}^2}{\left(s_N+\sigma^2\right)^2}\left \Vert \sum_{j=\ell+1}^{\ell_{GP}} \lambda_j \pvec{j} \pvec{j}^T  + \sigma^2 \mat{I} \right \Vert_2^2  + \sum_{j=\ell+1}^{\ell_{GP}} \lambda_j^2 \pvec{j}^T \mat{R}^2 \pvec{j}  \right). \nonumber
\dedit
\end{equation}
\end{theorem}
\smallskip
\sedit 
We can consider a relative notion of error by dividing each side by $\Vert \mbf{y} \Vert_2^2$. Furthermore, since we are dealing with bounded functions $f$ and finite amounts of data, we can always normalize and center the data to make $\Vert \mbf{y} \Vert_2 = \mathcal{O}(1)$.
\dedit
The difference between the two approximations described in Theorem~\ref{thm:GPS} has two sources. First is the contribution of kernel eigenfunctions that are not in the pseudospectral approximation basis, i.e., the summations involving $j > \ell$ above. Second is the contribution of the noise term $\sigma^2$.

Because the kernel eigenfunctions span the RKHS containing the GP posterior mean, the case of equivalent approximation spaces for pseudospectral approximation and GP regression follows simply by setting $\ell_{GP} = \ell < \infty$. \sedit This case is highlighted by Corollary~\ref{cor:equiv}.\dedit

\smallskip
\begin{corollary}{\label{cor:equiv}}
    \sedit
    Let $\ell_{GP} = \ell$. Then the difference between the pseudospectral approximation $\hat{f}_{\ell}(x)$ and  GP posterior mean $m(x)$ defined in Theorem~\ref{thm:GPS} is:
    \begin{equation}
        \Vert m - \hat{f}_{\ell}\Vert_{L^2_{\mu}}^2 \leq \Vert \mbf{y} \Vert_2^2 \, \ell N M^2 w^2_{\max}\left( \frac{\sigma^2}{s_N+\sigma^2}\right)^2 .
    \end{equation}
    \dedit
\end{corollary}
\smallskip
\sedit
In this case, the difference between approximations is due only to the nugget $\varn > 0$ and the smallest eigenvalue of the covariance matrix as indicated by $s_N$. If we additionally have that $N \leq \ell_{GP}$, then the design covariance matrix $\mat{R}^{-1}$ remains invertible (i.e., with $s_N > 0$) even as $\varn \to 0$, yielding zero difference between the approximations. 
\dedit
This occurs, for example, in the case of kernels constructed from fully tensorized polynomial eigenfunctions and tensorized Gaussian quadrature rules, where $N= \ell = \ell_{GP}$.

We also note that the bound in Theorem~\ref{thm:GPS} depends on the minimum eigenvalue of the covariance matrix, 
 represented (after diagonalization) by $\mat{S} + \sigma^2\mat{I}$. If $\mat{S}$ is nearly rank-deficient and the nugget is sufficiently small, then the bound can be large. This situation is not purely an artifact of the theory; indeed, it corresponds to a poorly conditioned numerical problem, where the actual difference between the two approximations may be large as well. One may imagine a case where $\mat{R}$ is not invertible, e.g., too many quadrature nodes are used and the GP kernel is finite rank; in this case the computation of $m$ becomes unstable whereas the computation of $\hat{f}_{\ell}$ can still proceed in a stable manner.
%
%
%


Besides bounding the difference between approximations, we can also analyze the impact of an orthogonalizing experimental design on the integrated variance of the GP posterior. 
Consider, again, a kernel consisting of $\ell_{GP}$ eigenfunctions and a training set corresponding to the nodes of a quadrature rule $\mathcal{Q}$ that orthogonalizes the first $\ell$ eigenfunctions. We begin by splitting~\eqref{eq:intVarBasis} into summations involving the first $\ell$ eigenfunctions and the remaining $\ell+1$ to $\ell_{GP}$ eigenfunctions:
\begin{align*}
\int c(x) d\mu (x) &= \sum_{i=1}^{\ell}\lambda_i \left(1 - \lambda_i\pvec{i}^T  \mat{R} \pvec{i}\right) +  
            \sum_{i=\ell+1}^{\ell_{GP}}\lambda_i \left(1 - \lambda_i\pvec{i}^T  \mat{R} \pvec{i}\right) .
\end{align*}
The second term in this expansion represents the contribution of the extra eigenfunctions to the integrated variance. We cannot comment on how our training points will affect this term because they are designed to have special properties only for the first $\ell$ eigenfunctions. But we can use these properties to analyze the first term above, thus describing the impact of an orthogonalizing rule on the associated integrated variance. Note that the weights $\mat{W}$ do not explicitly enter the GP regression; nonetheless, as we show in Appendix~\ref{sec:appVar}, this portion of the integrated variance can be rewritten as:
\begin{align}
\label{eq:psaonivar}
\sum_{i=1}^{\ell}\lambda_i \left(1 - \lambda_i\pvec{i}^T  \mat{R} \pvec{i}\right) &= \sum_{i=1}^{\ell}\lambda_i \pvec{i}^T\mat{W}\left( \sum_{j=\ell+1}^{\ell_{GP}}\lambda_j\pvec{j}\pvec{j}^T + \sigma^2 \mat{I} \right) \mat{R} \pvec{i}.
\end{align}
This expression consists of interactions between the first $\ell$ eigenfunctions and the last eigenfunctions; it also includes the impact of the nugget. The eigenfunction interactions are expected because the design only orthogonalizes the first $\ell$ eigenfunctions; errors are incurred when the numerical inner product is taken between one of the first $\ell$ eigenfunctions and one of the remaining eigenfunctions. We see that if $\ell=\ell_{GP}$, the integrated variance is of the order of the noise $\sigma^2$ as expected. When additionally $\sigma^2 = 0$, the integrated variance is zero. An orthogonalizing design ensures that the integrated variance captures only the contributions of eigenfunctions which are not orthogonalized.

The preceding results highlight some of the assumptions underlying the practical application of these two surrogate construction methodologies. When using a pseudospectral approximation, one implicitly assumes that the true function's projection onto basis functions \textit{not} included in the expansion (and hence not orthogonalized by the underlying design) is small. Otherwise, more points and basis functions should be added to the approximation; indeed, adaptive basis selection is the concern of a vast array of approximation methods, pseudospectral and otherwise. 
In GP regression, a properly chosen kernel is one whose eigenvalue decay \sedit matches the decay of the spectral coefficients of the true function. \dedit In this case the function will lie in the RKHS associated with the kernel and an accurate approximation can readily be achieved. 
These two approaches interact through the integrated variance. The orthogonalizing rule in a pseudospectral approximation ensures that the difference between $f$ and $f_{\ell}$ lies mostly in the span of the eigenfunctions $\phi_{k>\ell}$. Correspondingly, this rule forces the IVAR, a measure of the uncertainty and error (via the Bayes risk) of the GP approximation, to retain contributions only from the same subspace (that spanned by $\{ \phi_{k>\ell} \}$). We also see that the difference between the GP posterior mean and  pseudospectral approximation is dominated by these extra eigenfunctions.

\subsection{Relaxing exact orthogonality}
\label{sec:spectra}
It is instructive to consider the tradeoff between exactly orthogonalizing \textit{fewer} basis functions or, for the same design points, generating an approximation using more basis functions than can be orthogonalized. The latter corresponds to performing GP regression under the conditions of Theorem~\ref{thm:GPS} with $\ell_{GP} > \ell$. 

Consider the function $f(x) = \sin(\pi x + 0.2)$ for $x \in \mathbb{R}$ with standard Gaussian weight $\mu$. We will use a Gauss-Hermite quadrature rule with $N=20$ points to construct a pseudospectral approximation of $f$ with the first $\ell=20$ Hermite polynomials (degrees 0 to 19) as basis functions. We compare this approximation to Gaussian process regression on the same points using the Mehler kernel~\cite{Schaback2007}
\begin{equation}
    K(x, y;t) = \frac{1}{\sqrt{1-t^2}} \exp\left(-\frac{1}{2}\frac{\left(x\right)^2t^2 - 2txy + \left(y\right)^2t^2}{1-t^2}\right),
\label{eqn:mehler}
\end{equation}
which is a closed-form expression for $K(x,y) =  \sum_{i=0}^{\infty}t^i\phi_i(x) \phi_i(y)$, where $\{\phi_i \} $ are again normalized Hermite polynomials. The nugget  $\varn = 0$ and the decay constant $t = 0.8$. In this way we compare two surrogates using identical experiments and closely related approximation spaces: the approximation space of the former is contained in that of the latter. But the two surrogates use different approaches for combining the same function evaluations into an approximation. GP regression approximates $f$ using an infinite collection of polynomial eigenfunctions, with more emphasis on those associated with higher eigenvalues, while pseudospectral approximation projects $f$ onto finite polynomial basis, with the exactness of the projection limited by external aliasing.

\begin{figure}[ht!]
\center
\includegraphics[scale=1.0]{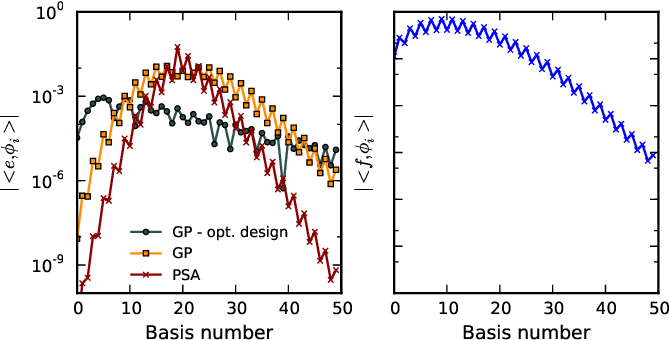}
\caption{Spectrum of approximation errors for the example function discussed in Section~\ref{sec:spectra}. Left figure shows the projection of the approximation error onto normalized Hermite polynomial basis functions $\phi_i$, for three different surrogates: pseudospectral approximation with 20 Gauss quadrature points, GP regression with a Mehler kernel on the same Gauss quadrature points, and GP regression with a Mehler kernel on IVAR-optimal points. Right figure shows the magnitude of the exact projection of $f$ onto each basis function $\phi_i$.}
\label{fig:psavsgpproj}
\end{figure}

Figure~\ref{fig:psavsgpproj} shows the results of this experiment. First, we note that the relative $L^2_{\mu}$ errors are $8.7 \times 10^{-3}$  for pseudospectral approximation and $1.8 \times 10^{-3}$ for the GP mean function. Perhaps more interesting than this improvement, however, is the \textit{spectrum} of the error associated with each approximation. The left panel of Figure~\ref{fig:psavsgpproj} shows the projections of the errors $e_{ps}(x) = \hat{f}_{\ell}(x)-f(x)$ and $e_{gp}(x) = m(x) -f(x)$ onto each basis function $\phi_i$, i.e., $|\langle e, \phi_i \rangle |$. For reference, the right panel shows how much energy $f$ has in each basis direction, via the magnitudes of the \textit{exact} projections $|\langle f, \phi_i \rangle |$. (All projections are computed with extremely high-order quadrature.) We see that the projection of the pseudospectral approximation error $e_{ps}$  onto the first few basis functions is small, even though $f$ itself has significant energy in these directions. The projection of $e_{ps}$ then rises with the basis index and peaks around an index of 20; it then begins to decay, just as the exact coefficients $\langle f, \phi_i \rangle$ themselves decay in magnitude. The projection of the GP error, on the other hand, is flatter across the indices. It rises slightly for the first few basis functions and then begins to decay somewhat more slowly. The basis functions for which the orange line (GP error) lies below the red line (PSA error) correspond to relatively large coefficient values; this results in a lower $L^2$ error overall.

As a preview of the next section, we include a third line in Figure~\ref{fig:psavsgpproj} representing GP regression performed with a 20-point IVAR-minimizing design. The relative $L^2_\mu$ error of this approximation is even smaller: $1.0 \times 10^{-5}$. And the spectrum of its error, shown with the grey line in the left panel, is even flatter than that of the quadrature-based GP approximation. This result amplifies the previous trend: compared to pseudospectral approximation, GP regression spreads its error more evenly over the spectrum. Departing from an orthogonalizing design allows this error to be spread even more broadly.

%% file: GPPSAexp.tex
\subsection{Similar approximation spaces, optimal experiments}\label{sec:numCompare}

Our next goal is to compare GP regression and pseudospectral approximation when the approximation spaces are similar, but the experimental designs (and the algorithms) differ. We have already seen that RKHS containing the GP mean can coincide with the range of the pseudospectral approximation operator in the case of a finite-rank GP kernel: the eigenfunctions of the kernel can simply be used as the basis for the pseudospectral approximation. In the numerical comparisons below, however, we would like to allow the GP kernel to have infinite rank. This choice is more representative of how GP regression is used in practice, and in principle may allow the GP mean to better approximate a wider variety of functions. Of course, the eigenvalues of the kernel need to decay, so that we have a bounded kernel.

We will therefore use the Mehler kernel \eqref{eqn:mehler} as in Section~\ref{sec:spectra}, and consider target functions whose inputs are endowed with standard Gaussian weight on $\mathbb{R}^d$. The eigenfunctions of the Mehler kernel are Hermite polynomials, which we use as basis functions for pseudospectral approximation. In $d>1$ dimensions, we use a tensorized version of the Mehler kernel: $K(x,y;t_1, \ldots, t_d) = \prod_{i=1}^d K (x^{(i)},y^{(i)}; t_i )$,
where $t_i$ now governs the decay rate of the eigenvalues associated with the univariate eigenfunctions in dimension $i$.

\sedit
To design experiments for GP regression, we use IVAR minimization in combination with adaptation of the covariance kernel hyperparameters (e.g., correlation lengths) and the nugget $\sigma^2$. In particular, we choose these parameters by maximizing the log-marginal likelihood,
\begin{equation}
    \log p(\mbf{\hat{y}}| \mbf{x}, \theta) = -\frac{1}{2} \mbf{\hat{y}}^T\mat{R}\mbf{\hat{y}}  - \frac{1}{2}\log|\mat{R}^{-1}| - \frac{N}{2}\log 2\pi,
\end{equation}
with respect to the complete set of hyperparameters (including $\sigma^2$), denoted by $\theta$. Our approach interleaves kernel adaptation with batch IVAR design. The batch size for each design is described in each example, but the number of Monte Carlo points used to evaluate the IVAR objective is fixed at $N_{mc} = 10^4$. 
\dedit

For pseudospectral approximation, we use a state-of-the-art adaptive Smolyak algorithm~\cite{patrick}, which adaptively enriches both the approximation basis and the experimental design using a greedy heuristic. The approximation is essentially a Smolyak sum of full tensor-product polynomial approximations, each computed using a pseudospectral approach. We use this approach because a regular tensor-product approach becomes infeasible for more than a few dimensions; hence many sparse grid~\cite{Barthelmann2000} and dimension-adaptive~\cite{Gerstner2003} polynomial approximation algorithms have been developed. The Smolyak algorithm in~\cite{patrick} uses generalized sparse grids and dimension adaptivity and thus provides a useful benchmark for comparison with kernel-adaptive GP. It is implemented in the MIT Uncertainty Quantification framework (MUQ) \cite{muq}. Since the input domain is endowed with Gaussian measure, the sparse grid design is based on one-dimensional Gauss-Hermite quadrature rules, with the number of points growing exponentially with the level of the sparse grid.

\subsubsection{Additively separable functions}

\begin{figure}[h!]
\centering
\subfigure[Relative error comparisons]{
\includegraphics[scale=1.0]{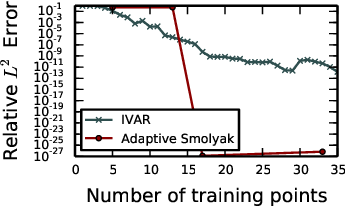}
\label{fig:easyPSAErr}
}
\subfigure[Trace of hyperparameter estimates]{
\includegraphics[scale=1.0]{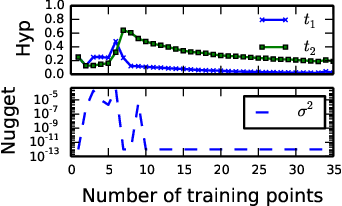}
\label{fig:easyPSATrace}
}
\caption{Error comparisons between GP/IVAR and adaptive Smolyak approximations for target function (\ref{eq:psa1}). Right panel shows hyperparameter traces for the GP design procedure.}
\label{fig:easyPSA}
\end{figure}

\begin{figure}
\centering
\subfigure[Relative error comparisons]{
\includegraphics[scale=1.0]{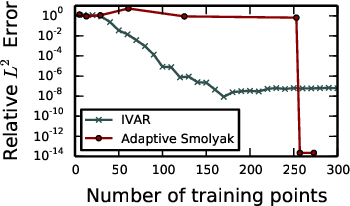}
\label{fig:sinPSAErr}
}
\subfigure[Trace of hyperparameter estimates]{
\includegraphics[scale=1.0]{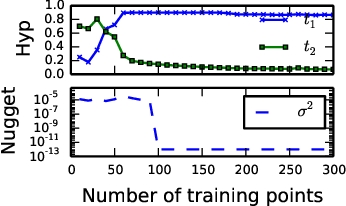}
\label{fig:sinTrace}
}
\caption{Error comparisons between GP/IVAR and adaptive Smolyak approximations for target function (\ref{eq:psa2}). Right panel shows hyperparameter traces for the GP design procedure.}
\label{fig:sinPSA}
\end{figure}

The first problems we consider are two dimensional and additively separable. Additive separability is a property favorable to the Smolyak algorithm, in part because of the greedy heuristic the algorithm uses to enrich the polynomial basis. The target functions are:
\begin{align}\label{eq:psa1}
f_1(x) &= x^{(1)} + (x^{(2)})^2, \\
f_2(x) &= \sin(2\pi x^{(1)}) + (x^{(2)})^2 \label{eq:psa2}.
\end{align}
The batch size for the closed-loop IVAR scheme is one training point. Adaptation is performed for the eigenvalue decay rates $t_1$ and $t_2$ in each dimension and for the noise term $\sigma^2$. Figures~\ref{fig:easyPSA} and~\ref{fig:sinPSA} show the resulting relative errors and hyperparameter traces, as a function of the number of training points. As in Section~\ref{sec:numerics}, the relative errors are defined as $||f - \hat{f}_{\ell}||_{L^2_{\mu}} / ||f||_{L^2_{\mu}}$ and $||f - m(x) ||_{L^2_{\mu}} / ||f||_{L^2_{\mu}}$ for Smolyak pseudospectral approximation and GP/IVAR, respectively. These errors are computed using 10000 Monte Carlo points.

In these examples, we observe that the adaptive Smolyak algorithm requires several function evaluations until it begins to converge. Once it converges, however, the relative error drops to machine precision. (Recall that the quadrature rules used in the Smolyak scheme require adding several points at a time.) The GP approximation error, on the other hand, decreases steadily but gradually as additional points are added and as the kernel is refined. In both cases, by the time the pseudospectral approximation becomes accurate, the GP approximation has already achieved at least $10^{-4}$ relative error, perhaps sufficient for many applications. The hyperparameter traces show how the GP covariance kernel adapts to the function being approximated. Function $f_1$ \eqref{eq:psa1} is fairly low order in both dimensions, and the hyperparameter values $t_1$ and $t_2$ both converge to fairly low numbers, indicating fast eigenvalue decay. The converged value of $t_1$ is small because only the constant and linear eigenfunctions (in the first dimension) are updated by the data; $t_2$ converges to a slightly higher value, thus slowing the eigenvalue decay, to account for the quadratic term in $x^{(2)}$. Function $f_2$ \eqref{eq:psa2}, on the other hand, is relatively high order in the first dimension but low order in the second. We thus observe that $t_1$ converges to $0.9$, corresponding to a slower decay, and that $t_2$ converges to roughly $0.1$; the hyperparameter optimization has ``learned'' that only the first few eigenfunctions in this dimension matter.

\subsubsection{Non-additively separable function}
\begin{figure}
\centering
\subfigure[Relative error comparisons]{
\includegraphics[scale=1.0]{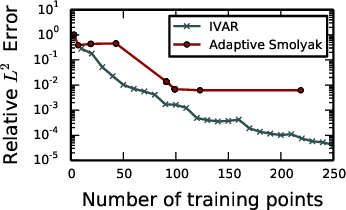}
}
\subfigure[Trace of hyperparameter estimates]{
\includegraphics[scale=1.0]{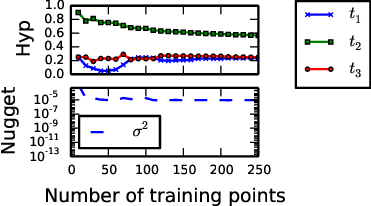}
}
\caption{Error comparisons between GP/IVAR and adaptive Smolyak approximations for the Ishigami function \eqref{eq:ishigami1}.}
\label{fig:ishigami}
\end{figure}

Next we perform the same experiments on the three-dimensional Ishigami function 
\begin{align} \label{eq:ishigami1}
    f(x) = \sin(x^{(1)}) + 7\sin^2(x^{(2)}) + 0.05\left(x^{(3)}\right)^4\sin(x^{(1)}),
\end{align}
The Ishigami function is not additively separable, and we expect to see even better relative performance of GP regression in this example since the kernel eigenfunctions include the tensor products of all univariate Hermite polynomials. Here we use a batch size of $M=10$ experiments for IVAR design and hyperparameter adaptivity. The results are shown in Figure~\ref{fig:ishigami}.

We see that the GP begins to outperform the pseudospectral approximation after roughly 20 function evaluations, with an error that decreases consistently. The error from the adaptive Smolyak approach, on the other hand, plateaus in several regions. After 50 iterations, the GP approximation reaches a relative error of $10^{-2}$, while it takes pseudospectral approximation 100 iterations to reach the same relative error. This behavior may be attributed to the fact that the Mehler kernel accounts for interactions among the input variables immediately, whereas the dimension-adaptive Smolyak approach requires some exploration (corresponding to the error plateaus) to find the basis functions that capture these interactions. 

\subsubsection{Ten-dimensional Genz function}
\begin{figure}
\center
\includegraphics[scale=1.0]{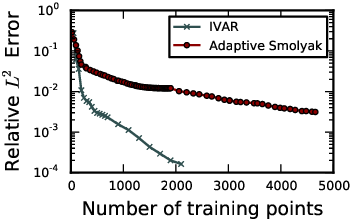}
\caption{Error comparisons between GP/IVAR and adaptive Smolyak approximations for the ten-dimensional Genz function~(\ref{eq:10Dgenz}).}
\label{fig:10Dgenz}
\end{figure}

Finally, we consider a higher-dimensional example using the oscillatory Genz function \cite{genz}:
\begin{equation} \label{eq:10Dgenz}
f_4(x) = \cos\left(2\pi w_1 + \sum_{i=1}^{10}x^{(i)} c_i\right),
\end{equation}
where $w_1=0.3$ and the $c_i$ are chosen randomly from a uniform distribution on $[0,1]$ and then normalized to $\Vert c \Vert_1 = 2.25$. This Genz function is typically evaluated on a hypercube domain $[-1,1]^{10}$ with normalization $\Vert c \Vert_1=9$. We have found the present approximation problem, with an unbounded and Gaussian-weighted domain, to be significantly more challenging, however, due to the tail behavior of the Hermite polynomials. 
A comparison of approximation errors, again between GP/IVAR using the Mehler kernel and adaptive Smolyak pseudospectral approximation, is shown in Figure~\ref{fig:10Dgenz}. The GP approximation performs extremely well. We note that this Genz function involves non-additive coupling among all ten input dimensions, a feature that may amplify the benefits of including a fully tensorized set of eigenfunctions via the GP kernel. For the GP regression calculations, we initially added experiments in batches of $M=50$, learning the hyperparameters after each batch, until obtaining 700 experiments total. Then we added experiments $200$ at a time. In the future, it may be useful to automatically stop adapting the hyperparameters once the changes in the hyperparameters between iterations fall below some threshold.

Note that we were able to compare GP and pseudospectral approximation in these examples because we had access to the Mehler kernel, and hence explicit evaluations of the Hermite basis functions were not required for GP regression. In practice, we may not have closed-form kernels for other basis function families, typically used in pseudospectral approximation. {Truncated} kernels whose eigenfunctions are polynomials, however, can often be represented using the Christoffel-Darboux formula~\cite{Christoffel1858, Darboux1878}. For normalized basis functions, the Christoffel-Darboux formula is essentially a finite rank kernel with equal eigenvalues, and thus for a fixed basis order, one would not be able to adapt this kernel.

%% file: conclusion.tex
\section{Conclusion}

This paper has examined experimental design criteria and optimization procedures used to develop surrogates of computational models, highlighting aspects of the interplay between experimental designs and approximation algorithms.
In the first part of the paper, we discussed experimental design criteria for Gaussian process regression and presented an algorithm for minimizing the posterior integrated variance (IVAR) of a Gaussian process over a continuous design space. 
\sedit
Our approach is adapted to the kernel structure; for instance, with isotropic squared exponential kernels, it yields well-spaced points on arbitrary complex domains, avoiding boundary clustering and other undesirable artifacts. 
\dedit
IVAR points also have good interpolation stability, as measured by the Lebesgue constant for kernel interpolation. The underlying optimization problem is tractably solved with gradient descent methods, as long as the number of design points is not too small or too large relative to the complexity of the prior covariance kernel. Our numerical experiments also demonstrate that simultaneously designing multiple points can yield substantial benefits over greedy strategies. Nonetheless, even greedy minimization of IVAR yields better approximation performance than standard greedy algorithms for minimizing entropy or maximizing MI.

In the second part of the paper, we compare GP regression to polynomial approximation. For simplicity, we consider only pseudospectral approximation, using nodes and weights that orthogonalize the finite set of basis functions used for approximation. In this setting, when GP regression uses the same nodes and a kernel with polynomial eigenfunctions, the difference between the two approximations is due only to the GP ``nugget'' and truncation of the kernel. When instead coupled with an IVAR design, GP regression for an infinite-rank kernel sacrifices numerical orthogonality of any given set of eigenfunctions; compared to the pseudospectral approach, projection errors may be larger for the first few eigenfunctions but are more nearly constant over the eigenspectrum. This observation is reminiscent of the average-case quadrature of \cite{Minka2000}, which compares a Bayesian method for deriving quadrature nodes with Gauss rules of fixed degree. We follow these comparisons with an empirical study of approximation performance on various test functions, comparing adaptive variants of GP regression and sparse pseudospectral approximation. We observe that while additively separable functions lend themselves easily to sparse polynomial approximations, more strongly coupled functions can be more efficiently approximated using the GP approach. Kernel approximation is also well suited to complex (e.g., non-tensorized) input domains. In our current design approach, while eigenfunctions of the kernel integral operator on the domain need not be explicitly computed, eigenfunctions that couple the input dimensions are all implicitly present.

Future work can extend our analysis of experimental design for GP regression in many ways. It is useful to compare the present design criteria with other methods for choosing ``good'' interpolation nodes in the radial basis function literature \cite{de2003optimal,de2005near}. Comparing IVAR-optimal nodes to optimal nodes for ``Bayesian quadrature'' \cite{Minka2000,Osborne2012} would also be of great interest. More broadly, we advocate for closer connections between the numerical analysis literature and the statistics literature on these topics. Finally, adaptive designs that interleave point selection with updates to the kernel would benefit from more rigorous study. Current sequential approaches, including the ones presented here, are perhaps reasonable heuristics. But look-ahead strategies that anticipate information gain from future designs, and that balance information for kernel adaptation with reduction of the conditional posterior variance, may lead to even more effective approximation procedures.

%% file: proofsContinuousVersion.tex
\newcommand{\dplus}{\pmat{\ell+1}{\ell_{GP}}\Lambda_{\ell+1:\ell_{GP}}\pmat{\ell+1}{\ell_{GP}}^T}
\newcommand{\dreg}{\pmat{2}{\ell}\Lambda_{1:\ell}\pmat{1}{\ell}^T}
\newcommand{\pw}{\pmat{1}{\ell}^T\mat{W}}
\newcommand{\wpt}{\mat{W}\pmat{1}{\ell}}
\newcommand{\si}{\mat{U}\left(\mat{S}+\sigma^2\mat{I}\right)^{-1}\mat{U}^T}
\newcommand{\sii}{\mat{U}\left(\mat{S}+\sigma^2\mat{I}\right)^{-2}\mat{U}^T}
\newcommand{\pwwp}{\pmat{1}{\ell}^T\mat{W}^2\pmat{1}{\ell}}
\newcommand{\pl}{\pmat{1}{\ell}\Lambda_{1:\ell}}
\newcommand{\lp}{\Lambda_{1:\ell} \pmat{1}{\ell}^T}
\newcommand{\px}[1]{\phi_{#1}\left(\mbf{x} \right)}
\newcommand{\ddp}{\sum_{i=\ell+1}^{\ell_{GP}}\lambda_i \pvec{i} \pvec{i}^T}
\newcommand{\ddpn}{\sum_{i=\ell+1}^{\ell_{GP}}\lambda_i \Vert\pvec{i} \pvec{i}^T\Vert^2}
\newcommand{\sumphif}{\sum_{j=1}^{\ell_{GP}}\lambda_j \pvec{j}} %
\newcommand{\sumphip}{\sum_{j=1}^{\ell}\lambda_j \pvec{j}} %
\newcommand{\sumphie}{\sum_{j=\ell+1}^{\ell_{GP}}\lambda_j \pvec{j}} %

\section{Proof of Theorem~\ref{thm:GPS}}\label{sec:errApp}
First we recall some notation. The prior precision matrix is $\mat{R} = \left(\sum_{j=1}^{\ell_{GP}}\lambda_j\pvec{j}\pvec{j}^T + \sigma^2 \mat{I}\right)^{-1}$. We define $\mat{U}$ and $\mat{S}$ to be the matrices associated with the eigenvalue decomposition $\mat{R} = \mat{U}\left(\mat{S} + \sigma^2 \mat{I}\right)^{-1}\mat{U}^T$.

We assume, without loss of generality, that the prior mean is zero; from~\eqref{eq:postMean} we obtain
\begin{align*}
m(x) = \mbf{y}^T \mat{R} \sumphif \phi_j (x) & = \mbf{y}^T \mat{R} \sumphip \phi_j (x) + \mbf{y}^T \mat{R} \sumphie \phi_j (x) \\
     &= \mbf{y}^T \mat{R} \sumphip \phi_j (x)  + a(x),
\end{align*}
where we have replaced $K(\mbf{x},x)$ with its eigenfunctions and then separated the posterior mean into two terms, the first containing the first $\ell$ eigenfunctions and the second $a(x)$ denoting the rest.

We now use the orthogonality rule~\eqref{eq:orthoQuad} to obtain $m(x) = \mbf{y}^T \mat{R} \sumphip \pvec{j}^T\mat{W}\pvec{j} \phi_j (x) + a(x)$. Recall that the pseudospectral approximation $\hat{f}_{\ell}$ is given by $\hat{f}_{\ell}(x) = \mbf{y}^T \mat{W} \sum_{j=1}^\ell\pvec{j} \phi_j(x)$, such that the difference between the GP and the spectral expansion becomes:
\begin{align*}
m(x) - \hat{f}_{\ell}(x) = \sum_{j=1}^{\ell} \left[ \mbf{y}^T\left( \mat{R} \lambda_j \pvec{j}\pvec{j}^T - \mathbf{I}\right) \mathbf{W}\pvec{j} \phi_j(x) \right] + a(x) .
\end{align*}
Now define $\mbf{d}^T := \mbf{y}^T\left(\mat{R}\sumphip \pvec{j}^T - \mat{I}\right) \mathbf{W}$ for convenience, and
form the $L_{\mu}^2$ difference
\begin{align}
\int \left(m - \hat{f}_{\ell}\right)^2 d\mu &=  \int \left(
                                              \sum_{j=1}^{\ell}
                                              \mbf{d}^T \pvec{j}
                                              \phi_j(x) + a(x)
                                              \right)^2 d\mu
                                              \nonumber \\
	&= \left \langle \sum_{j=1}^{\ell} \mbf{d}^T \pvec{j}
          \phi_j(x), \sum_{j=1}^{\ell} \mbf{d}^T \pvec{j} \phi_j(x)
          \right \rangle + \left \langle \sum_{j=1}^{\ell} \mbf{d}^T
          \pvec{j} \phi_j(x), a(x)  \right \rangle + \langle a(x),
          a(x) \rangle \nonumber \\
	&= \left \langle \sum_{j=1}^{\ell} \mbf{d}^T \pvec{j}
          \phi_j(x), \sum_{j=1}^{\ell} \mbf{d}^T \pvec{j} \phi_j(x)
          \right \rangle  +  \langle a(x), a(x) \rangle \label{eq:diffdecomp}
\end{align}
where the the cross terms resulting from the expansion of the square
disappear due to orthogonality between the first $\ell$ basis
functions and the $ k> \ell$ basis functions.  

We now seek to bound the size of the first term in \eqref{eq:diffdecomp}. Use the orthogonality properties of $\phi_j$ to obtain
\begin{align*}
\left \langle \sum_{j=1}^{\ell} \mbf{d}^T \pvec{j} \phi_j(x), \sum_{j=1}^{\ell} \mbf{d}^T \pvec{j} \phi_j(x) \right \rangle  &= \sum_{i,j=1}^{\ell} \pvec{i}^T \mbf{d} \mbf{d}^T\pvec{j} \langle \phi_i, \phi_j \rangle = \sum_{i=1}^{\ell}  \left(\mbf{d}^T\pvec{i} \right)^2 
\end{align*}
We can further simplify this expression by taking advantage of the invertibility of $\mat{R}$ and rewriting $\mbf{d}^T$ as
\begin{align*}
\sum_{i=1}^{\ell}  \left(\mbf{d}^T\pvec{j} \right)^2  &= \sum_{i=1}^{\ell} \left( \mbf{y}^T  \mat{R} \left( \sumphip \pvec{j}^T - \mat{R}^{-1} \right) \mat{W} \pvec{i} \right)^2 \\
 &=  \sum_{i=1}^{\ell} \left( \mbf{y}^T \mat{R} \left( \sumphie \pvec{j}^T - \sigma^2 \mat{I} \right) \mat{W} \pvec{i} \right)^2  \\
 &\leq \mbf{y}^T \mbf{y} \sum_{i=1}^{\ell}  \left \Vert \mat{R} \left( \sumphie \pvec{j}^T - \sigma^2 \mat{I} \right) \mat{W} \pvec{i} \right \Vert_2^2,
\end{align*}
where the first equality is obtained by extracting $\mat{R}$, the second equality results from substituting in the definition of $\mat{R}^{-1}$, and the bound is due to the Cauchy-Schwarz inequality. Now we can factorize the last term above because all induced norms are sub-multiplicative,
\begin{align*}
\sum_{i=1}^{\ell}  \left(\mbf{d}^T\pvec{j} \right)^2  &\leq \mbf{y}^T \mbf{y} \left \Vert \mat{R} \right \Vert^2 \sum_{i=1}^{\ell} \left \Vert  \left( \sumphie \pvec{j}^T  - \sigma^2 \mat{I} \right) \mat{W} \pvec{i} \right \Vert_2^2 \\
 &\leq  \mbf{y}^T \mbf{y} \left  \Vert\left(\mat{S}+\sigma^2\mat{I}\right)^{-1} \right \Vert^2 \left \Vert \sum_{j=\ell+1}^{\ell_{GP}} \lambda_j  \pvec{j} \pvec{j}^T  + \sigma^2\mat{I}\right \Vert_2^2 \sum_{i=1}^{\ell} \Vert\mat{W}\pvec{i}\Vert^2 \\
\sedit
 &\leq   \mbf{y}^T \mbf{y}  \frac{1}{\left(s_N + \sigma^2\right)^2}  \left \Vert \sum_{j=\ell+1}^{\ell_{GP}} \lambda_j  \pvec{j} \pvec{j}^T  + \sigma^2\mat{I}\right \Vert_2^2 \ell w^2_{\max} \pvec{z_2}^T\pvec{z_2} 
 \dedit
\end{align*}
\sedit where for the third equality we use the fact that $\mat{W}$ is
a diagonal positive definite matrix; $w^2_{\max}$ is the square of the
largest-magnitude weight,
$w^2_{\max} = \arg\max_{i \in \{1,\ldots,N\}}w^2_i$; and $z_2$ is the
index of the discrete basis function with maximum norm,
$z_2 = \arg \max_{i \in \{1, \ldots, \ell\}} \pvec{i}^T\pvec{i}$. Let
$M = \underset{(x,w) \in \mathcal{Q}_{\ell}, i \in \{1,\ldots
  \ell\}}{\max}|\phi_i(x)|$,
since the first $\ell$ eigenfunctions are bounded at the $N$
quadrature nodes. Then we have $ \pvec{z_2}^T\pvec{z_2} \leq M^2 N$
and
\begin{equation*}
\sum_{i=1}^{\ell}  \left(\mbf{d}^T\pvec{j} \right)^2 \leq   \mbf{y}^T \mbf{y}  \frac{\ell N M^2 w^2_{\max}}{\left(s_N + \sigma^2\right)^2}  \left \Vert \sum_{j=\ell+1}^{\ell_{GP}} \lambda_j  \pvec{j} \pvec{j}^T  + \sigma^2\mat{I}\right \Vert_2^2 
\end{equation*}
\dedit

We now turn to the $\langle a(x), a(x) \rangle$ term. This term
is simplified similarly to the first, using Cauchy-Schwarz and the
orthogonality of the basis functions:
\begin{align*}
\langle a(x), a(x) \rangle &\leq  \mbf{y}^T \mbf{y} \sum_{j=\ell+1}^{\ell_{GP}}\lambda_j^2 \pvec{j}^T  \mat{R} \mat{R} \pvec{j} .
\end{align*}

Thus the overall bound can be given as:
\begin{align*}
    \sedit
    \Vert m - \hat{f}_{\ell}\Vert_{L^2_{\mu}}^2 \leq \mbf{y}^T \mbf{y} \left(\frac{\ell N M^2 w_{\max}^2}{(s_N+\sigma^2)^2}  \left \Vert \sum_{j={\ell}+1}^{\ell_{GP} }\lambda_j \pvec{j} \pvec{j}^T  + \sigma^2 \mat{I} \right \Vert_2^2  + \sum_{j=\ell+1}^{\ell_{GP}} \lambda_j^2 \pvec{j}^T \mat{R}^2 \pvec{j}  \right).
    \dedit
\end{align*}

\section{Derivation of Equation~(\ref{eq:psaonivar})}\label{sec:appVar}

We begin by splitting \eqref{eq:intVarBasis} into two components in order to focus on the first term:
\begin{align}\label{eq:apintcovsplittwo}
\int c d \mu &= \sum_{i=1}^{\ell}\lambda_i \left(1 - \lambda_i\pvec{i}^T  \mat{R} \pvec{i}\right) + \sum_{i=\ell+1}^{\ell_{GP}}\lambda_i \left(1 - \lambda_i\pvec{i}^T  \mat{R}  \pvec{i}\right).
\end{align}

We again use the trick of replacing one with $\pvec{i}^T\mat{W} \pvec{i}$ to obtain
{ \small
\begin{align*}
\sum_{i=1}^{\ell}\lambda_i \left(1 - \lambda_i\pvec{i}^T  \mat{R}  \pvec{i}\right) &= 
    \sum_{i=1}^{\ell}\lambda_i \left(\pvec{i}^T\mat{W}\pvec{i} - \lambda_i \pvec{i}^T  \mat{R}  \pvec{i} \right) \\
    &= \sum_{i=1}^{\ell}\lambda_i \pvec{i}^T\left(\mat{W}- \lambda_i  \mat{R}  \right) \pvec{i} \\
    &= \sum_{i=1}^{\ell}\lambda_i \pvec{i}^T\left(\mat{W}\left( \sum_{j=1}^{\ell_{GP}}\lambda_j \pvec{j} \pvec{j}^T +\sigma^2\mat{I}\right)- \lambda_i \mat{I} \right) \mat{R} \pvec{i} \\
    &= \sum_{i=1}^{\ell}\lambda_i \pvec{i}^T\left( \lambda_i \mat{I} + \mat{W} \sum_{j=\ell+1}^{\ell_{GP}}\lambda_j\pvec{j}\pvec{j}^T + \sigma^2 \mat{W} - \lambda_i \mat{I} \right) \mat{R} \pvec{i} \\
    &= \sum_{i=1}^{\ell}\lambda_i \pvec{i}^T\mat{W}\left( \sum_{j=\ell+1}^{\ell_{GP}}\lambda_j\pvec{j}\pvec{j}^T + \sigma^2 \mat{I} \right) \mat{R} \pvec{i},
\end{align*} }
where the first equality comes from orthogonality, the second equality results from extracting the basis vectors, the third equality results from extracting $\mat{R}$ on the right, and the fourth equality comes from splitting the sum over all $\ell_{GP}$ terms into two sums: from $1$ to $\ell$ and from $\ell+1$ to $\ell_{GP}$.

%% file: gradient_appendix.tex
\section{Gradient-based IVAR minimization}
\label{sec:gradientappendix}
\sedit The gradient of the IVAR objective \eqref{eq:mcobjective} with respect to the design variables can be obtained analytically from the underlying covariance kernel; this enables a wide variety of optimization approaches to be used for IVAR minimization. For simplicity, assume that we are designing experiments on a one-dimensional input domain and that our experimental design is $\mbf{x} = [x_1,\ldots x_N]$. Now the gradient for a given $\hat{x}_i$ (see \eqref{eq:mcobjective}) may be computed as  
\begin{align}
  \nabla_{\mbf{x}} c(\hat{x}_i | \mbf{x}) &= -\nabla_{\mbf{x}}\left( K(\mbf{x},\hat{x}_i)^T\mat{R}K(\mbf{x},\hat{x}_i)\right) \\
    &= -\nabla_{\mbf{x}}\left(\sum_{j=1}^N \sum_{k=1}^N K(x_j,\hat{x}_i)\mat{R}[j,k]K(x_k,\hat{x}_i)\right) \nonumber
\end{align}
Let us focus on a single element of the gradient:
\begin{align*}
  \frac{\partial c(\hat{x}_i | \mbf{x})}{\partial x_{\ell}} &= -\left[ \sum_{j=1}^N \sum_{k=1}^N \frac{\partial}{\partial x_{\ell}} \left( K(x_j,\hat{x}_i) \mat{R}[j,k] K(x_k,\hat{x}_i)\right) \right] \\
&= - \sum_{j=1}^N \sum_{k=1}^N \left( \frac{\partial}{\partial x_{\ell}}  K\left(x_j,\hat{x}_i\right) \right) \mat{R}[j,k] K(x_k,\hat{x}_i) - \sum_{j=1}^N\sum_{k=1}^N   K\left(x_j,\hat{x}_i\right) \frac{\partial}{\partial x_{\ell}}\left( \mat{R}[j,k] K(x_k,\hat{x}_i\right) ) \\
 &= 2 \left(\frac{\partial}{\partial x_{\ell}} K(x_{\ell}, \hat{x}_i)\right)\sum_{k=1}^N\mat{R}[\ell,k]K(x_k,\hat{x}_i) - \sum_{j=1}^N\sum_{k=1}^N   K\left(x_j,\hat{x}_i\right) \frac{\partial}{\partial x_{\ell}}\left( \mat{R}[j,k]\right)  K(x_k,\hat{x}_i) 
\end{align*}
where the second equality follows from the chain rule and the third equality follows from the symmetry of $K$ and $\mat{R}$ as well as another application of the chain rule. We are left with terms involving derivatives of $K(x,x^{\prime})$ and of the elements of $\mat{R}$. To evaluate the latter, use the identity $\frac{\partial \mat{R}}{\partial x_{\ell}} = -\mat{R} \frac{\partial \mat{R}^{-1}}{\partial x_{\ell}} \mat{R}$ and recall that $\mat{R}^{-1}[i,j] = K(x_i,x_j) + \delta_{ij} \sigma^2 $. Hence this quantity can also be computed from the derivative of $K.$ For instance, if $K$ is a squared exponential kernel $K(x_i,x_j) = \exp\left (-(x_i-x_j)^2/2l^2 \right )$ then the derivative is ${\partial K(x_i,x_j)}/{\partial x_i} = -\frac{1}{l^2}(x_i-x_j)K(x_i,x_j).$
An analogous derivation can be performed for a multi-dimensional input space.
